\RequirePackage[leqno]{amsmath}
\documentclass[leqno]{amsart}
\usepackage{amsmath,mathtools}
\makeatletter
\usepackage[utf8]{inputenc}     
\usepackage{multicol}
\usepackage[inline]{enumitem} 
\usepackage{color, graphics}
\usepackage[demo]{graphicx}
\usepackage{efbox,graphicx}
\efboxsetup{linecolor=gray,linewidth=0.8pt}
\usepackage{floatrow}
\usepackage{float}
\usepackage{pgfplots}
\usepackage{mathrsfs}
\usepackage{latexsym, amsmath, amssymb}
\usepackage{pgfplots}
\pgfplotsset{compat=newest}
\usepackage{graphicx}
\usepackage{subfigure}
\usetikzlibrary{backgrounds,automata}
\usepackage{caption}
\usepackage{pgfplots}
\usepackage{amsfonts}
\usepackage{pifont}
\usepackage{enumitem}
\setlist[itemize]{topsep=0pt,after=\vspace{1.5\baselineskip}}
\usepackage[colorlinks=true]{hyperref}
\usepackage{empheq}
\usepackage[T1]{fontenc}
\usepackage{tabularx}
\usepackage[leftcaption]{sidecap}
\sidecaptionvpos{figure}{m}
\usepackage{tikz}
\usepackage{caption}
\usepackage{tikz-cd}
\usetikzlibrary{decorations.markings}
\usetikzlibrary{decorations.pathmorphing}
\usepackage{tabularx,multirow,array} 
\usepackage[margin=0.5in]{geometry}
\usetikzlibrary{backgrounds,automata}

\pgfplotsset{compat=1.18}

\usepackage{xparse}
\ExplSyntaxOn
\int_new:N \g_tohi_const_int
\int_new:N \g_tohi_const_sub_int
\tl_new:N  \g_tohi_const_char_tl
\cs_new_protected:Nn \tohi_print_constant:nn 
{
  #1 \textsubscript {#2} 
} 

\NewDocumentCommand\resetconstants{m}
{
 \int_gincr:N \g_tohi_const_int
 \int_gzero:N \g_tohi_const_sub_int
 \tl_gset:Nn  \g_tohi_const_char_tl {#1}
}

\NewDocumentCommand\const{m}
{
  \tl_if_exist:cTF
   {
    c_tohi_const_\int_use:N\g_tohi_const_int _#1_tl
   }
   {
    \tl_use:c {c_tohi_const_\int_use:N\g_tohi_const_int _#1_tl }
   }
   {
    \int_gincr:N \g_tohi_const_sub_int
    \tl_const:cx {c_tohi_const_\int_use:N\g_tohi_const_int _#1_tl }
     { \exp_not:N\tohi_print_constant:nn {\g_tohi_const_char_tl }{\int_use:N \g_tohi_const_sub_int}}
    \tl_use:c {c_tohi_const_\int_use:N\g_tohi_const_int _#1_tl }
   }
}

\ExplSyntaxOff
\newcommand{\inlineitem}[1][]{%
\ifnum\enit@type=\tw@
    {\descriptionlabel{#1}}
  \hspace{\labelsep}%
\else
  \ifnum\enit@type=\z@
       \refstepcounter{\@listctr}\fi
    \quad\@itemlabel\hspace{\labelsep}%
\fi}
\DeclarePairedDelimiter\abs{\lvert}{\rvert}

\DeclarePairedDelimiter\tonda{(}{)}
\DeclarePairedDelimiter\quadra{[}{]}

\newcommand{\into}{\int_\Omega}

\setlist[itemize]{noitemsep, topsep=0pt}

\def\R{\mathbb R} \def\N{\mathbb N}

\def\R{\mathbb R} \def\N{\mathbb N} 
\def\TM{T_{\textup{max}}} 
\def
\@cite
#1#2{[#1\if@tempswa, #2\fi]}
\makeatother
\newtheorem{theorem}{Theorem}[section]

\newtheorem{lemma}[theorem]{Lemma}

\newtheorem{remark}{Remark}

\title[Boundedness and blow-up for Keller--Segel models in penetrable domains] 
{
To what extent does the consideration of positive total flux influence the dynamics of Keller--Segel-type models?
}
\author[Khadijeh Baghaei, Silvia Frassu, Yuya Tanaka, Giuseppe Viglialoro]{$^{\dagger}$Khadijeh Baghaei,
$^{\sharp}$Silvia Frassu, $^{\flat}$Yuya Tanaka, and $^{\natural}$Giuseppe Viglialoro$^{\star}$}
\makeatletter\@namedef{subjclassname@2020}{\textup{2020} Mathematics Subject Classification}\makeatother
\keywords{Chemotaxis, Robin-type boundary conditions, Incoming flow, Gradient nonlinearities, Boundedness. \\
\textit{$^\star$Corresponding author}:giuseppe.viglialoro@unica.it}
\begin{document}
\maketitle
{
\centerline{$^{\dagger}$Pasargad Institute for Advanced Innovative Solutions}
\centerline{No. 30, Hakim Azam St., North Shiraz St., Mollasadra Ave., Tehran, Iran}
\medskip
\centerline{$^{\flat}$Department of Mathematical Sciences}
\centerline{Kwansei Gakuin University}
\centerline{1, Uegahara, Gakuen, Sanda, Hyogo, 669--1330, Japan}
\medskip
\centerline{$^{\natural}$Dipartimento di Matematica e Informatica}
\centerline{Universit\`{a} degli Studi di Cagliari}
\centerline{Via Ospedale 72, 09124. Cagliari, Italy}
\medskip
}
\bigskip
\resetconstants{c}
\begin{abstract}
Since the introduction of the Keller-Segel model in 1970 to describe chemotaxis (the interactions between cell distributions $u$ and chemical distributions $v$), there has been a significant proliferation of research articles exploring various extensions and modifications of this model within the scientific community. From a technical standpoint, the totality of results concerning these variants are characterized by the assumption that the total flux, involving both distributions, of the model under consideration is \textit{zero}. This research aims to present a novel perspective by focusing on models with a \textit{positive} total flux. Specifically, by employing Robin-type boundary conditions for $u$ and $v$, we seek to gain insights into the interactions between cells and their environment, uncovering important dynamics such as how variations in boundary conditions influence chemotactic behavior. In particular, the choice of the boundary conditions is motivated by real-world phenomena and by the fact that the related analysis reveals some interesting properties of the system. 

Mathematically, or $h,\chi>0$ and $\alpha \in (0,1]$ we investigate Keller--Segel-type models with positive total flux $u_{\nu}-\chi u v_{\nu}=\chi \alpha h uv$, reading as  
\begin{equation}\label{problemAbstract}
\tag{$\oplus$}
\begin{dcases}
u_t= \Delta u -\chi \nabla \cdot (u \nabla v) & \textrm{in } \Omega \times (0,\TM),\\
\tau v_t= \Delta v - v  +u & \textrm{in } \Omega \times (0,\TM),\\
u_{\nu}=(\alpha-1)\chi h u v, \; v_{\nu}=-h v & \textrm{on } \partial \Omega \times (0,\TM),\\
\end{dcases}
\end{equation}
where $\Omega \subset \mathbb{R}^n$, $n \geq 1$ is a bounded and smooth domain, $\tau\in \{0,1\}$, $\TM>0$ and $\nu$ denoting the outward normal vector to the boundary of $\Omega$, $\partial \Omega$.

We aim at emphasizing how the inclusion of the incoming flowing flow makes the overall analysis more complex.
The technical difficulties are essentially tied to the lack of the crucial property of the mass conservation, which in this case is replaced by an increase in the mass itself. Such behavior of the mass cannot be circumvented by merely including classical logistics of the form  $a u-b u^2$ (with $a,b>0$); an additional dissipative term involving gradient nonlinearities is required. But there is another indication that suggests how the total positive flux significantly alters the dynamics of taxis models with zero-flux. Indeed: 
%
\begin{itemize}
\item [$\lozenge$] for mechanisms with positive flux ($\alpha>0$), global (i.e. $\TM=\infty$) and bounded solutions are obtained as long as an essential strong logistic of the form $a u-b u^2-c|\nabla u|^2$ ($a,b,c>0$) is included in the model; 
%
\item [$\lozenge$] for phenomena with zero flux ($\alpha=0$), analogies with the classical Keller--Segel models, without logistic, endowed with homogeneous Neumann boundary conditions, specially for blow-up scenarios (i.e., $\TM<\infty$), can be observed. 
\end{itemize}
\end{abstract}
\tableofcontents
\section{Chemotaxis phenomena. The homogeneous Neumann boundary conditions}\label{SectionTheLandamrkingKS}
\subsection{The landmarking Keller--Segel model}
The term \textit{chemotaxis} indicates the movement in a certain environment of living organisms (e.g. cells) occupying habitats under the influence of a chemical stimulus, therein distributed and changing over time. 
The pioneering mathematical model for the description of the mentioned phenomenon was proposed by physicist Keller and mathematician Segel in 1970 (\cite{K-S-1970,Keller-1971-MC}), by means of two coupled Partial Differential Equations. If $(x,t)$ identifies a position $x$ in a region (a bounded domain $\Omega$ of $\R^n$, with $n\geq 1$) at some time $t$ (with $t\in (0,\TM)$, and $\TM\in (0,\infty]$), the classical Keller--Segel model reads:
\begin{equation}\label{KS}
\begin{cases}
 u_t = \Delta u -\chi \nabla \cdot(u\nabla v), \qquad v_t=\Delta v-v+u  & \textrm{ in } \Omega \times \left(0,\TM\right),\\
u_\nu=v_\nu=0\;\textrm{ on } \;\partial \Omega \times \left(0,\TM\right), \;
u(x,0)=u_0(x), \;  v(x,0)= v_0(x)& x \in \bar\Omega.
\end{cases} 
\end{equation}
This system describes the self-organized aggregation of a starving slime mould colony. In this context, individual cells (with density $u=u(x,t)\geq 0$ and initial configuration $u_0(x)=u(x,0)\geq 0$) direct their otherwise random movement toward higher concentrations of an attractive chemical signal substance, which they produce (with concentration $v=v(x,t)\geq 0$ and initial distribution $v_0(x)=v(x,0)\geq 0$). 
Here, $\partial \Omega$ represents the boundary of $\Omega$, and $\nu$ is its outward normal vector. The homogeneous Neumann boundary conditions (also referred to as ``h.N.b.c.'') $u_\nu=v_\nu=0$ on $\partial \Omega \times (0,\TM)$ idealize an insulated and impenetrable domain. 

The cell density equation in system \eqref{KS} reflects the balance between diffusion ($\Delta u$) and the amplifying effect of taxis-driven chemosensitivity ($-\chi\nabla \cdot (u\nabla v)$), which intensifies with larger $\chi$, due to the absence of external sources and the isolation of the domain, any net growth in cell mass over time is excluded. This follows from integrating over $\Omega$ and applying the h.N.b.c., yielding $\frac{d}{dt} \int_\Omega u = 0$, so $\int_\Omega u$ stays equal to its initial value $\int_\Omega u_0(x)dx$.

Mathematically, model \eqref{KS} admits both global bounded solutions (with $\TM = \infty$ and $\lVert u(\cdot,t)\rVert_{L^\infty(\Omega)}$ finite for all $t > 0$) and unbounded solutions that blow up in finite time ($\TM$ finite, with $\limsup_{t\to \TM} \lVert u(\cdot,t)\rVert_{L^\infty(\Omega)} = \infty$). For $n=1$, all solutions remain bounded over time (diffusion dominates self-attraction), while for $n \ge 2$, self-attraction can overpower diffusion, leading to finite-time blow-up depending on the initial data. Classical results and surveys, \textit{all derived by relying on the mentioned boundedness of the mass}, are available, for instance, in \cite{BellomoEtAl,HerreroVelazquez,Nagai,WinklAggre,LankeitWinkler-FacingLow2020}.
\subsection{The Keller--Segel system with logistic and strong logistic sources}\label{SectionNilTotalFlux}
Especially in a setting involved with population dynamics, the combination of \eqref{KS} with sources distributed \textit{internally} in the domain describing population growth or decay is very natural; we are specially referring to logistic terms (see \cite{AiTsEdYaMi,OsTsujYagMimuBidim}). 

In this context, species distribution is also affected by reproduction and mortality, with death rates increasing at higher population densities. The general \textit{logistic source} in the equation
$
u_t = \nabla u - \chi \nabla \cdot (u \nabla v) + f(u)
$
takes the form
$
f(u) = a u^{\alpha} - b u^{\beta}
$, where \( a, b \geq 0 \) and \( 1 \leq \alpha < \beta \); the positive term models births, while the negative term models deaths. We now examine various scenarios:
\begin{itemize}
\item [$\triangleright$]
First, for $b=0$ and $a>0$, a spatial integration of the equation for the cells implies, due to the h.N.b.c (and using the H\"{o}lder inequality) 
\begin{equation}\label{BagadheiRelazBlow}
\frac{d}{dt}\int_\Omega u= a \int_\Omega u^\alpha\geq a |\Omega|^{1-\alpha}\left(\int_\Omega u\right)^\alpha \quad \textrm{for all } t\in (0,\TM),
\end{equation} 
which leads to explosion of $\int_\Omega u$ (and henceforth of $\lVert u(\cdot,t)\rVert_{L^\infty(\Omega)}$) at 
\[
\TM\leq \frac{|\Omega|^{\alpha-1} \left(\int_\Omega u_0(x)dx\right)^{1-\alpha}}{a (\alpha-1)}.
\]
\item [$\triangleright$] Secondly, for $a=0$ and $b>0$, similarly, we have 
\[
\frac{d}{dt}\int_\Omega u= -b \int_\Omega u^\beta \leq 0 \quad \textrm{for all } t\in (0,\TM),
\] 
thus, the mass decreases; $\displaystyle \int_\Omega u\leq \int_\Omega u_0(x)dx$ for all $t\in (0,\TM)$.
\item [$\triangleright$] For $a,b>0$ and $\beta>\alpha\geq1$ the uniform-in-time finiteness of $\sup_{t\in (0,\TM)}\int_\Omega u$ can be as well achieved with some rearrangements, entailing
\[
\int_\Omega u \leq \max\left\{\int_\Omega u_0(x)dx, \left(\frac{a}{b}|\Omega|^{\beta-\alpha}\right)^\frac{1}{\beta-\alpha}\right\} \quad \textrm{for all } t \in (0,\TM).
\]
More precisely, analogously to the previous cases, we obtain 
\[
\frac{d}{dt} \int_\Omega u= a \int_\Omega u^\alpha - b \int_\Omega u^\beta \quad \textrm{for all } t \in (0,\TM),
\]
and, thanks to the H\"{o}lder inequality, we have
\[
-\into u^\beta \leq -\abs*{\Omega}^{\frac{\alpha-\beta}{\alpha}}\tonda*{\into u^\alpha}^\frac{\beta}{\alpha} \quad \textrm{and} \quad -\left(\into u^\alpha\right)^\frac{\beta-\alpha}{\alpha} \leq -\abs*{\Omega}^\frac{(1-\alpha)(\beta-\alpha)}{\alpha}\tonda*{\into u}^{\beta-\alpha}\quad \text{in }(0,\TM).
\]
By combining the above expressions, we arrive at 
\[
\frac{d}{dt} \int_\Omega u \leq a \int_\Omega u^\alpha \left(1 - \frac{b}{a} \abs*{\Omega}^{\alpha-\beta} \left(\int_\Omega u \right)^{\beta-\alpha}\right) \quad \textrm{on }(0,\TM),
\]
so concluding by invoking the ODE comparison principle in \cite[Lemma 3.3]{ChiyoDuzgunFrassuVigliaoro-2024}.
\end{itemize}
A natural question is whether the inclusion of logistic terms can ensure boundedness and prevent finite-time blow-up detected in \eqref{KS}. This has been confirmed only for sufficiently large \( b \) (e.g., \( \beta = 2 \), see \cite{TelloWinkParEl, W0}), while blow-up may still occur for certain \( \beta > 1 \) close to 1 (see \cite{WinDespiteLogistic, Winkler_ZAMP-FiniteTimeLowDimension}). However, blow-up can be avoided by strengthening the logistic source. Specifically, if \( f(u) \) is replaced by \( h(u,|\nabla u|) = f(u) - c|\nabla u|^\gamma \), then any solution to
$
u_t = \Delta u - \chi \nabla \cdot (u \nabla v) + h(u,|\nabla u|)
$
remains globally bounded for all \( c > 0 \), \( \beta > \alpha \geq 1 \), provided \( \frac{2n}{n+1} < \gamma \leq 2 \) (see \cite{IshidaLankeitVigliloro-Gradient, LiEtAl2024gradientnonlinearitiesprevent} for theoretical and applied results on gradient-dependent sources in chemotaxis).  
\begin{remark}[On the logistic terms with strong/gradient dissipative effects]\label{RemarkStrongLogistic}
The mention of the logistic term $h=h(u,|\nabla u|) $with gradient nonlinearities at this stage of the discussion is intentional and not without purpose. Indeed, these terms will prove to be of crucial importance in our forthcoming discussion, as will become apparent in $\S$\ref{SectionLogistiStrongLogistic} below.
\end{remark}
\subsection{Nonlinear diffusion and sensitivity and further models} 
Over recent decades, numerous extensions of model \eqref{KS} have been proposed to describe taxis-driven phenomena, often formulated as
\begin{equation}\label{KS-General}
 u_t = \nabla \cdot (S(u,v)\nabla u - T(u,v)\nabla v) + f(u), \qquad v_t = \Delta v + g(u,v) \quad \text{in } \Omega \times (0,\TM),
\end{equation}
with homogeneous Neumann boundary conditions and given initial data. Here, \( S(u,v) \) models diffusion, \( T(u,v) \) captures chemotactic sensitivity (positive for attraction, negative for repulsion), and \( f(u) \) accounts for source effects. The chemical \( v \) may be produced (\( g(u,v) = -v + u \)) or consumed (\( g(u,v) = -uv \)) by the cells.

Depending on $S,T,f,g$, solutions may exhibit either global boundedness or finite-time blow-up. For parameters $m_1, m_2\in \R$, $\rho,\gamma>0$ and $\beta > \alpha \geq 1$, we briefly review key studies, the associated qualitative behaviors of the aforementioned terms:
\begin{itemize}
\item [$\triangleright$] Nonlinear diffusion and sensitivity model with generalized logistic and nonlinear production: $S(u,v)\simeq u^{m_1}, T(u,v)\simeq u^{m_2}, f(u)\simeq u^\alpha-u^\beta, g(u,v)\simeq -v+u^\gamma$ (\cite{BellomoEtAl,Lankeit,verhulst,AiTsEdYaMi,TanakaViglialoroYokota}). 
\item [$\triangleright$] Linear consumption model:
$S(u,v)= 1, T(u,v)= u, f(u)\equiv 0, g(u,v)= -vu$
 (\cite{TaoBoun,BaghaeiaKhelghatib}). 
\item [$\triangleright$] Chemotaxis model with flux limitation and linear production:
$S(u,v)=u/\sqrt{u^2+|\nabla u|^2}, T(u,v)=u/\sqrt{1+|\nabla v|^2}, f(u)\equiv 0, 
g(u,v)=-v+u$ or $g(u,v)=-\frac{1}{|\Omega|} \int_\Omega u + u$
(\cite{Bellomo04032017, BellomoWinklerFlux, ChiyodaEtAlFlux,MarrasVernierYokota,MIZUKAMI20195115, WinklerFlux2022}). 
\item [$\triangleright$] Chemotaxis model with singular (or logarithmic) sensitivity and nonlinear consumption:
 $S(u,v)=1, T(u,v)=u/v, f(u)\equiv 0, g(u,v)\simeq -u^\gamma v$
 (\cite{LankeitViglialoroAAM,WinklerRenorm,WinklerSingSens}).
\item [$\triangleright$] Linear production model with local sensing diffusion: 
$S(u,v)\simeq v^{-\rho}u, T(u,v)=u, f(u)\equiv 0, g(u,v)=-v+u $
(\cite{FujieSenbaSensing,FujieCalVar}). 
\item [$\triangleright$] Nonlinear diffusion and sensitivity model with nonlocal source  and linear production:
$S(u,v)\simeq u^{m_1},T(u,v)\simeq u^{m_2}, f(u)=u^\alpha(1-\int_\Omega u^\beta), g(u,v)=-v+u$ (\cite{ChiyoDuzgunFrassuVigliaoro-2024,NegreanuEtAl}).
\end{itemize}
\textit{To clarify our presentation, it is crucial to note that the results discussed so far depend on the key property of cellular distribution: the uniform-in-time boundedness of mass, which is closely linked to the homogeneous Neumann boundary conditions.}
\section{Modification in the boundary conditions. Requiring the total flux to vanish}
Recently, research on chemotactic systems has expanded to explore various types of boundary conditions, moving beyond the Neumann-type conditions analyzed in previous sections. \textit{Our project is situated within this evolving framework}; in particular, the foundations of our investigation are based on the the following considerations.
\subsection{The uniform-in-time $L^1(\Omega)$-bound of $u$} \label{BiundenessMass-Section}
As noted in $\S$\ref{SectionTheLandamrkingKS}, the key property for solutions of chemotactic problems is the uniform-in-time $L^1(\Omega)$-boundedness. This is immediately satisfied under homogeneous Neumann boundary conditions (\( u_\nu = v_\nu = 0 \) on \( \partial \Omega \times (0, \TM) \)), or when zero \textit{total flux} is imposed under other boundary conditions. Specifically, for model \eqref{KS-General}, the total flux is given by the vector field \( S(u,v)\nabla u - T(u,v)\nabla v \), and zero total flux is enforced by the boundary condition
\begin{equation}\label{ZeroTotalFluxBoundary}
S(u,v) u_\nu - T(u,v) v_\nu = 0 \quad \text{on} \quad \partial \Omega \times (0, \TM).
\end{equation}
Table \ref{Tabella} summarizes key results from the literature on chemotaxis mechanisms, highlighting specific expressions for \( S \), \( T \), \( h \), and \( g \), alongside various boundary conditions and zero total flux. As anticipated, the condition \eqref{ZeroTotalFluxBoundary} ensures that integrating the first equation of model \eqref{KS-General} over \( \Omega \) avoids boundary integral terms. Specifically, we have $
\frac{d}{dt} \int_\Omega u = \int_{\partial \Omega} \left( S(u,v) u_\nu - T(u,v) v_\nu \right) + \int_\Omega f(u) = \int_\Omega f(u)$ on $(0,\TM)$, which implies the necessary uniform-in-time boundedness of \( \lVert u(\cdot,t) \rVert_{L^1(\Omega)} \).
\begin{table}[h!]
\begin{tabular}{m{0.78cm}||m{1.18cm} m{1.18cm} m{1.18cm} m{1.18cm}||m{3.98cm}}
& $S(u,v)$ & $T(u,v)$ & $f(u)$ &$g(u,v)$ & Boundary condition for $v$\\ 
 \hline
 (1) & $1$ & $u$ & &$u$& $v_{\vert \partial \Omega\times (0,\TM)}=0$ \\
(2)&$1$ & $u$& & $-v+u$ &  $v_{\vert \partial \Omega\times (0,\TM)}=0$  \\
(3)&$1$ & $u$& & $-vu$ &  $v_{\vert \partial \Omega\times (0,\TM)}=v^*$  \\
(4)&$1$ & $u$& $u-u^2$& $-vu$ &  $v_{\nu\vert \partial \Omega\times (0,\TM)}=1-v$  \\
(5)&$1$ & $u$& $u-u^2$& $-vu$ &  $v_{\vert \partial \Omega}=v^*$  \\
(6)&$u^m$ & $u$& & $-vu$ &  $v_{\nu\vert \partial \Omega\times (0,\TM)}=1-v$  \\
(7)&$u^m$ & $-u/v$& & $-vu$ &  $v_{\vert\partial \Omega\times (0,\TM)}=v^*$  \\
(8)&$u^m$ & $-u$& & $-vu$ &  $v_{\vert\partial \Omega\times (0,\TM)}=v^*$  \\
\end{tabular}
\caption{(1) Analysis of a two-dimensional system with Dirichlet boundary conditions for which solutions remain bounded near the boundary and the blowup set has a finite number of interior points (\cite{SuzukiDirichlet2013}); 
(2) Study of critical mass levels in $n$-dimensional domains (\cite{FuhrmannLankeitWinklerDirichlet});
(3) Asymptotics and stationary solution for chemotaxis-consumption models (\cite{BraukhoffLankeit2019,YangAhn2024});
(4)   Existence of a unique bounded classical solution for a chemotaxis-Navier--Stokes model involving an incompressible fluid. (In the table we avoid the interaction with fluid.) (\cite{BraukhoffMarcelNS-2017}); 
(5) Long time behavior of solutions to chemotaxis-consumption systems with logistic growth (\cite{knosallalankeit2024}); 
(6)   Asymptotics for a nonlinear chemotaxis-Navier--Stokes model involving an incompressible fluid. (In the table we avoid the interaction with fluid.) (\cite{ChunyanZhaoyin2020});
(7) Existence of blowing-up solution at finite time (\cite{Wang-Winkler-2023});
(8) Detection of critical exponent $m$ in two-dimensional settings (\cite{AhnWinkler-2023}).
 }
 \label{Tabella}
\end{table}
 \section{Positive total flux in chemotaxis. Real applications and mathematical formulation}\label{Intro}
\subsection{The Robin boundary conditions}
In the literature, boundary conditions of the form \( r \psi + s \psi_\nu = g \) on \( \partial \Omega \), where \( \psi \) and \( g \) are functions and \( r, s \) are constants, are known as Robin (or mixed) boundary conditions (possibly abbreviated as ``R.b.c.''). For specific values of \( (r,s) = (1,0) \) or \( (r,s) = (0,1) \), these reduce to Dirichlet or Neumann conditions, respectively. In their more general form, Robin boundary conditions provide important advantages for modeling complex phenomena. For instance, as described by Souplet in \cite{Souplet_Gradient}, for a single (biological) species density \( u \) in a bounded domain, a Robin boundary condition models a scenario in which the species can cross the boundary of its environment. (Further examples and applications of Robin boundary conditions in parabolic problems can be found in \cite{AndreGia, SunYudong, QS, PP, PS}.)

In the frame of taxis-driven processes where the response of cells to chemical gradients near boundaries is a crucial aspect, Robin boundary conditions may offer a more nuanced representation of boundary interactions compared to Dirichlet or Neumann conditions. As illustrated in Table \ref{Tabella}, R.b.c. have already been employed in the study of chemotaxis mechanisms, specifically in modeling the behavior of chemical signals under zero-flux conditions.
\subsection{The Bacillus subtilis, Eukaryotic T-Cells, and Microbial Communities}\label{SectionRoleChemo}
This study investigates how chemical concentrations influence total flux across boundaries, with particular attention to the role of boundary conditions in shaping chemotactic responses. In the below biological transport mechanisms, where total flux is governed by both chemical gradients and cellular distribution, a net positive flux often emerges—capturing the hallmark of active cellular adaptation to environmental stimuli.

In line with \cite{TuvalEtAl}, Bacillus subtilis often thrives in thin fluid layers near oxygen contact lines, where oxygen uptake rates are singular and controlled by boundary concentrations. 

Chemotaxis is commonly measured by tracking the number of cells crossing a boundary per unit time, such as cells migrating through a chemoattractant-infused filter. In \cite{SzatmaryNossal}, eukaryotic cell migration is used as an indicator of chemotactic responsiveness.

In tumor biology, models couple tumor cell dynamics with chemokine concentration, driving T-cell migration toward the tumor via blood vessels, with the movement rate proportional to chemoattractant levels (see \cite{AlmeidaEtAl}).

Swarms and biofilms, microbial communities forming at the interface of solid substrates and strong fluid flows, are also influenced by boundary dynamics. These communities interact with flow-induced signal gradients, underscoring the importance of boundary conditions in microbial behavior (see \cite{WINKLE-OpenBoundaries}).
\subsection{The role of the chemoattractant and positive total flux. Mathematical interpretation and formulation}\label{SectionRolePositiveFlux}
We now present a consistent mathematical formulation of the phenomena discussed above. Since in chemotaxis mechanisms the cells direct their movement in the exact opposite direction to that of the gradient of the chemical signal, if the homogeneous Neumann boundary conditions $u_\nu=v_\nu=0$ on $\partial \Omega \times (0,\TM)$ are replaced by 
\begin{equation*}\label{RobinvNeumanu}
u_\nu=0\quad \textrm{and}\quad v_\nu=-h v \quad \textrm{on}\quad \partial \Omega \times (0,\TM),
\end{equation*}
indicating an outflow of the chemoattractant and the null flux of the cells' distribution, we derive that the total flux $u_\nu-uv_\nu=h u v$ is \textit{nonnegative} on $\partial \Omega \times (0,\TM)$. This effectively models how the outward flux of the chemoattractant \( v \) triggers cellular migration across the boundary into the domain, resulting in a \textit{inward total flux.}  

But the relevance of this transport mechanism can be more importantly intuited by this observation, which will make meaningful the definition of our boundary conditions. Indeed, even for appropriate outward flux of the cells' configuration, the taxis-driven effect of the outward flow of the chemoattractant can yet keep producing a positive total flux. Formally, for $\alpha\in (0,1]$ the following boundary conditions on $u$ and $v$
\begin{equation}\label{RobinvRobinu}
u_\nu=(\alpha-1)\chi h u v \quad \textrm{and}\quad v_\nu=-h v \quad \textrm{on}\quad \partial \Omega \times (0,\TM),
\end{equation}
both modeling outward flux for $u$ and $v$, turn to an essentially positive total flux:
\begin{equation}\label{RobinvRobinuPositiveNetFlux}
u_\nu-\chi u v_\nu = \alpha \chi h u v \quad  \textrm{on}\quad \partial \Omega \times (0,\TM).
\end{equation}
Consequently, the boundedness of the mass discussed in $\S$\ref{BiundenessMass-Section} is no more ensured. Indeed, by integrating the first equation in \eqref{KS} and taking into consideration \eqref{RobinvRobinuPositiveNetFlux} it is only possible to conclude that 
\begin{equation}\label{IncreasingMass}
\frac{d}{dt}\int_\Omega u =\alpha \chi h \int_{\partial \Omega} uv \geq 0 \quad \textrm{for all } t \in (0,\TM).
\end{equation}
This relation merely indicates that the total mass \( \int_\Omega u \) exhibits a nondecreasing behavior over the interval \( (0, T_{\mathrm{max}}) \), but it does not establish a direct connection between the evolution of the mass and a functional dependence on the mass itself. In particular, if one were able to derive from \eqref{IncreasingMass} an ordinary differential inequality analogous to that in \eqref{BagadheiRelazBlow}, then (by the same reasoning) finite-time blow-up would necessarily follow. However, such an inference is not supported by the available data, and it remains entirely plausible that the total mass may instead converge to a finite value as \( t \to T_{\mathrm{max}} \).
\begin{remark}[On the connection between conditions \eqref{RobinvRobinu} and \eqref{RobinvRobinuPositiveNetFlux}]
The parameter \( \alpha \) effectively regulates the balance between the fluxes of \( u \) and \( v \), as defined in equation \eqref{RobinvRobinu}, and their contribution to the total flux in equation \eqref{RobinvRobinuPositiveNetFlux}. When \( \alpha \to 0 \), the boundary permits substantial outflow of \( u \), and the advective influx due to \( v \) is insufficient to maintain a positive net flux. In contrast, for \( \alpha = 1 \), the outflow of \( u \) is completely suppressed, and the total flux reaches its maximum, driven solely by the outward transport induced by \( v \). Undoubtedly, the chosen boundary conditions reflect real-world phenomena and reveal valuable and interesting insights into the system's behavior. 
\end{remark}
\begin{figure}
\centering
\begin{tabular}[c]{ccc}
  \includegraphics[width=0.3\linewidth]{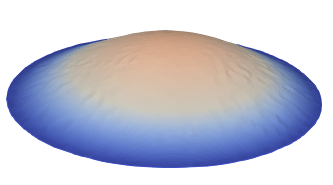}
  &\includegraphics[width=0.3\linewidth]{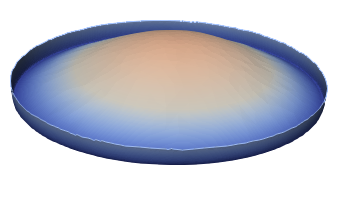}
  &\includegraphics[width=0.3\linewidth]{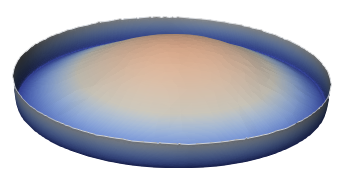}
  \\
  \includegraphics[width=0.3\linewidth]{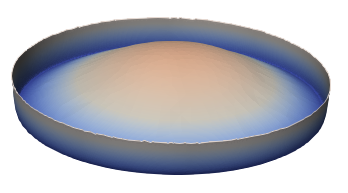}
  &\includegraphics[width=0.3\linewidth]{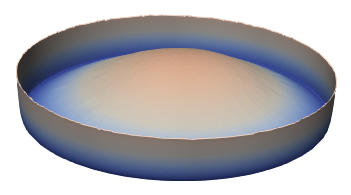}
  &\includegraphics[width=0.3\linewidth]{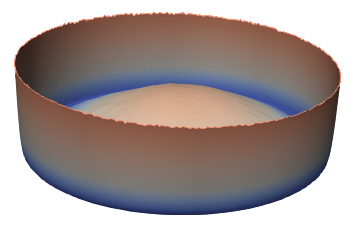}
\end{tabular}
  \caption{Simulations of model \eqref{KS} with boundary conditions given in \eqref{RobinvRobinu}. Evolution of the cells' density $u$ at six different instants of time: $t=0$ (top left), $t=0.1$ (bottom right). Some quantitative values are deductible from Figure \ref{FigureComparison}.}\label{fig:test-2d}
\end{figure}
\begin{figure}
\centering
\setlength{\subfigcapmargin}{5pt}
\subfigure[The solid (dotted) blue line represents the evolution of the  maximum of the solution for the Robin problem with $\alpha=1$ ($\alpha=0.7$), while the solid red line shows the solution for the Neumann problem. (Logarithmic scale is used to display the vertical axis.)]{
\begin{tikzpicture}\label{AFigure}
\begin{axis}[
width=9cm,height=8cm,
        xmin=0,xmax=550,
        ymin=0,ymax=700,
        ytick = {0,120,...,700},
        ymode=log, 
    xlabel=$t$,ylabel=$|| \cdot ||_{L^\infty(\Omega)}$]
  \addplot [blue,line width=0.35mm] table [x=a,y=b] {plotMaxR_L1.dat};
   \addplot [red,line width=0.35mm] table [x=a,y=b] {plotMaxN_L1.dat};
  \addplot [blue,dotted, line width=0.5mm] table [x=a,y=b] {plotMaxR_L07.dat};
\end{axis}
\end{tikzpicture}
}
\;\;\;\;
\subfigure[The solid (dotted) blue line represents the evolution of mass of the solution for the Robin problem with $\alpha=1$ ($\alpha=0.7$), while the solid red line shows the solution for the Neumann problem.]{
\begin{tikzpicture}\label{BFigure}
\begin{axis}[
width=9cm,height=8cm,
        xmin=0,xmax=550,
        ymin=10,ymax=70,
        ytick = {10,20,...,70},
    xlabel=$t$,ylabel=$|| \cdot ||_{L^1(\Omega)}$]
\addplot [blue,line width=0.35mm] table [x=a,y=b] {plotMassR_L1.dat};
  \addplot [red,line width=0.35mm] table [x=a,y=b] {plotMassN_L1.dat};
   \addplot [blue,dotted,line width=0.5mm] table [x=a,y=b] {plotMassR_L07.dat};
\end{axis}
\end{tikzpicture}
}
\caption{Robin conditions vs. Neumann conditions.}\label{FigureComparison}
\end{figure}
\begin{remark}[Some numerical evidences]
To support the theoretical insights discussed above, we present numerical experiments for the classical Keller--Segel system \eqref{KS} with Neumann boundary conditions and its variant with the modified boundary conditions in \eqref{RobinvRobinu}. The simulations are carried out on the unit disk,
\[
\Omega = \{(x, y) \in \mathbb{R}^2 : x^2 + y^2 < 1\},
\]
with initial conditions and parameters:
\[
u_0(x, y) = v_0(x, y) = 13 e^{-x^2 - y^2}, \quad \chi = 0.14, \quad h = 60.
\]
We focus on the case of positive total flux, corresponding to \( \alpha = 1 \) in \eqref{RobinvRobinuPositiveNetFlux}, where the boundary condition becomes \( u_\nu - \chi u v_\nu = \chi h u v \). Figure~\ref{fig:test-2d} shows how this setting leads to an accumulation of cells near the boundary over time. Figure~\ref{FigureComparison} compares the dynamics for different values of \( \alpha \). In Figure~\ref{AFigure}, the maximum of the cell density under the positive-flux model (\( u^P \)) is consistently greater than or equal to that of the zero-flux model (\( u^Z \)), with cusps marking moments when \( \max u^P > \max u^Z \). Figure~\ref{BFigure} compares the evolution of total mass in both cases, highlighting the impact of the boundary condition.
\end{remark}
\section{Boundedness of the mass by means of strong logistics in models with positive flux}\label{SectionLogistiStrongLogistic}
As previously discussed, in chemotaxis models with homogeneous Neumann boundary conditions, logistic terms of the form \( a u - b u^2 \) are commonly introduced to ensure the boundedness of \( u \) on the entire time interval \( (0, \infty) \). In contrast, such terms are not even sufficient to guarantee the boundedness of the total mass on \( (0, T_{\max}) \) when positive total flux boundary conditions are imposed. Indeed, in this case, relation~\eqref{IncreasingMass} takes the form
\begin{equation}\label{IncreasingMassBIS}
\frac{d}{dt} \int_\Omega u = \alpha \chi h \int_{\partial \Omega} u v + a \int_\Omega u - b \int_\Omega u^2, \quad \text{for all } t \in (0, T_{\max}),
\end{equation}
where the boundary integral \( \int_{\partial \Omega} u v \) remains difficult to control and may dominate the dynamics.

Given the critical role of uniform-in-time mass boundedness in the theoretical study of taxis-driven systems, it becomes essential to identify mechanisms that ensure this property under the boundary conditions introduced in \eqref{RobinvRobinu}. As we shall see, more robust and efficient damping terms are required to effectively counteract the potential growth of \( \int_\Omega u \) on \( (0, T_{\max}) \).
\subsection{Logistics with damping gradient terms; applications and their role through the boundedness of the mass}\label{SectionDampeningGrad}
As it will be shown, the mathematical estimates developed to control the evolution of the mass naturally inherit a structure that includes two key positive contributions: one involving the square of the cell density, and another involving the quadratic norm of its gradient. Consequently, to effectively manage the growth of \( \int_\Omega u \) described in \eqref{IncreasingMassBIS}, the introduction of suitably designed damping terms with gradient nonlinearities (as discussed in Remark~\ref{RemarkStrongLogistic}) becomes not merely a modeling choice, but a necessary component to ensure the qualitative stability and biological plausibility of the system under consideration.

Specifically, we will consider the modified source term \( h(u,|\nabla u|) = a u - b u^2 - c |\nabla u|^2 \) in the Keller--Segel model \eqref{KS}. Then, the cell density evolves according to the following equation:
\begin{equation*} \label{problem}
u_t = \Delta u - \chi \nabla \cdot (u \nabla v) + h(u, |\nabla u|) \quad \text{in } \Omega \times (0, T_{\max}).
\end{equation*}
At this stage, to ensure the boundedness of the mass, in this research we will operate in a twofold manner. The first approach relies on the use of trace inequalities, by means of which it will be possible to rewrite the undesired boundary integral as the sum of integrals over the entire domain. The second approach, which does not rely on the trace theory, involves the introduction of appropriate functionals and ad hoc manipulations, aimed at rendering any boundary terms negligible.

As will be rigorously established in the sequel, both strategies provide a mechanism for controlling the boundary integral corresponding to the positive part of the total flux. In essence, this quantity can be estimated in terms of a suitable norm of \( u \) in \( W^{1,2}(\Omega) \), at least morally:
\begin{equation*}
\alpha \chi h\int_{\partial \Omega} u v \leq \const{PD_1}\int_\Omega u^2+\const{PD_2}\int_\Omega |\nabla u|^2\quad \textrm{for all } t \in (0,\TM) \textrm{ and some } \const{PD_1}, \const{PD_2}>0.
\end{equation*} 
\begin{remark}[Quadratic-gradient-driven phenomena and their connection to positive fluxes]
While gradient nonlinearities arising from external sources, interpreted in biological contexts as mechanisms of ``accidental deaths'', have been previously studied within single-equation models (see, e.g., \cite{Souplet_Gradient,Souplet2005Book,QS,chipot_weissler,fila} and references therein), their impact within chemotactic systems has only recently begun to attract attention. Notably, recent works such as \cite{IshidaLankeitVigliloro-Gradient, LiEtAl2024gradientnonlinearitiesprevent} have started to explore their relevance in biologically and ecologically meaningful settings. In light of the preceding considerations (both mathematical and biological) it should no longer be surprising that, in models characterized by a positive total flux, there exists an intrinsic connection between the growth of the total mass and the presence of external agents exerting strong damping effects. 
\end{remark}
\section{Aim of the investigation. Models and results}
Motivated by both the biological insights and the mathematical considerations discussed above, and taking into account the apparent lack of comprehensive studies on chemotaxis models featuring \emph{positive total flux}, we set out to conduct a rigorous analysis of the existence, boundedness, and potential blow-up behavior of solutions to a class of chemotaxis systems. These systems are subject to Robin-type boundary conditions, as specified in equation~\eqref{RobinvRobinu}. The adoption of such boundary conditions introduces a broader and more flexible analytical framework, enabling us to capture a wider range of biological phenomena. Our investigation reveals the following key findings:

\begin{itemize}
    \item[$\lozenge$] In the case of mechanisms involving \emph{positive flux} (i.e., when $\alpha > 0$), we establish the global-in-time existence ($\TM = \infty$) and boundedness of solutions, provided that the system incorporates an essential strong logistic damping term of the form $a u - b u^2 - c|\nabla u|^2$ with positive constants $a, b, c > 0$. This damping mechanism plays a critical role in controlling solution growth and preventing blow-up. (Theorem \ref{theoremlocal}.)
    \item[$\lozenge$] When the flux is \emph{zero} (i.e., $\alpha = 0$), the model behavior draws significant analogies to the classical Keller--Segel models without logistic terms, which are typically studied under homogeneous Neumann boundary conditions. In such settings, we observe the emergence of blow-up scenarios (i.e., $\TM < \infty$), emphasizing the importance of flux and boundary effects in determining the qualitative behavior of solutions. (Theorems \ref{theoremlocalN1}, \ref{thm2} and \ref{BlowUpParabolic} in the final $\S$\ref{Appendix}.)
\end{itemize}
\subsection{Presentation of the models and of the main results}
Accordingly to the above preparations, we herein are interested to address questions connected to local-well posedness and global solvability of the following chemotaxis models with gradient-dependent sources and nonnegative total flux: 
\begin{equation}\label{problem}
\tag{$\mathcal{P}_\tau$}
\begin{dcases}
u_t= \Delta u -\chi \nabla \cdot (u \nabla v) + au-bu^2-c\abs{\nabla u}^2& \textrm{in } \Omega \times (0,\TM),\\
\tau v_t= \Delta v - v  +u & \textrm{in } \Omega \times (0,\TM),\\
u_{\nu}=(\alpha-1)\chi h u v, \; v_{\nu}=-h v & \textrm{on } \partial \Omega \times (0,\TM),\\
u(x,0)=u_0(x), \tau v(x,0)=\tau v_0(x)  & x \in \bar\Omega.
\end{dcases}
\end{equation}
Here, we assume 
\begin{equation}\label{reglocal}
    \textrm{ for some } \delta\in(0,1) \textrm{ and } n\in \N, \Omega \subset \R^n \textrm{ is a bounded domain of class } C^{2+\delta}.
\end{equation}
We will prove the following result. 
\begin{theorem}\label{theoremlocal}
Let the hypotheses in \eqref{reglocal} be fulfilled,  $\tau\in\{0,1\}$, $\chi, a, h >0$ and $\alpha \in [0,1]$. Then there exists $C_{\partial \Omega}=C_{\partial \Omega}(\Omega, n)$ such that for every
\begin{equation*}
 u_0, \tau v_0: \bar{\Omega}  \rightarrow \R^+ , \textrm{ with } u_0,  \tau v_0 \in C^{2+\delta}(\bar\Omega) \quad \textrm{complying with} \quad u_{0\nu}=(\alpha-1)\chi h u_0v_0 \quad \textrm{ and } \quad \tau v_{0\nu}=-h \tau v_0 \quad  \textrm{ on }\partial \Omega
\end{equation*}
the following conclusion holds true:
\begin{itemize} 
\item [] For $\tau=0$, whenever either 
\begin{equation}\label{Estimate_b}
\begin{cases}
b> c+ \frac{(\chi \alpha)^2 h C_{\partial \Omega}}{16 c}  & \text{ if } 0<c<\frac{\chi \alpha}{4} \sqrt{h C_{\partial \Omega}},\\
b>\frac{\chi \alpha}{2} \sqrt{h C_{\partial \Omega}} & \textrm{ if } c \geq \frac{\chi \alpha}{4} \sqrt{h C_{\partial \Omega}} \quad (\textrm{being } c>0 \,\textrm{ for }\, \alpha = 0),
\end{cases}
\end{equation}
or
\begin{equation}\label{Estimate_Bag0}
\alpha \chi <\min\{b,4c\}, 
\end{equation}
\item [] for $\tau=1$, whenever \eqref{Estimate_b} or
\begin{equation}\label{Estimate_Bag1}
b>\chi,
\end{equation}
\end{itemize}
 problem \eqref{problem} admits a unique solution
\begin{equation*}
(u,v)\in C^{2+\delta,1+\frac{\delta}{2}}( \Bar{\Omega} \times [0, \infty))\times C^{2+\delta,\tau+\frac{\delta}{2}}( \Bar{\Omega} \times [0, \infty))
\end{equation*}
such that $0\leq u,v \in L^\infty(\Omega \times (0,\infty)).$ 
\end{theorem}
Let us give some comments.
\begin{remark}
For $\alpha=0$, corresponding to the zero-flux situation, from relation \eqref{Estimate_b} it is seen that boundedness is ensured for $b, c>0$ only; despite that, we underline that the same property is achieved also for logistics with weaker powers in the dampening gradient terms. (See \cite[Theorem 1.2]{IshidaLankeitVigliloro-Gradient}.) Moreover, specifying $\alpha =0$ in the second assumption of \eqref{Estimate_b} is due to the fact that in such a limit case, the claim could not be ensured for $b > 0$ and $c = 0$ only. 
\end{remark}
\section{Necessary tools: parameters and properties of solutions to differential equations}
We will make use of the following results and we also underline that without any mention, all the constants $c_i$, $i=1,2,\ldots$ are assumed to be positive.

\subsection{Smoothing properties of Robin heat semigroup}
In view of analyses concerning regularity and estimates of solutions to parabolic problems with Robin conditions, we need to adapt known results for the same problems equipped, however, by homogeneous Neumann or Dirichlet conditions. In the specific, the following relations are called \textit{ultracontractivity} estimates. A characterization of semigroups that satisfy such estimates can be found in \cite[$\S$ 7.3.2]{Arendt-Wolfgang2004}. 
\begin{lemma}\label{LemmaRobinSemigroup}
Let $\Omega\subset\R^n$ satisfy condition in \eqref{reglocal} and $h>0$. For 
\[
D(\Delta_R) = \{ z \in H^1(\Omega) \mid \Delta z \in L^2(\Omega), \; z_\nu = -h z  \text{ on } \partial \Omega \},
\]
let us consider the semigroup of linear operators $(e^{t \Delta})_{t\geq 0}$ generated by the Robin Laplacian defined on $D(\Delta_R)$. Then there exist positive constants \( C_R \) and \( \mu_1 \) such that for all $1\leq q\leq p\leq \infty$ these inequalities hold for any $t>0$:
    \begin{equation}\label{SemigroupRobin1}
    \| e^{t\Delta} z \|_{L^{p}(\Omega)} \leq C_R \left(
1+t^{-\frac{n}{2}(\frac{1}{q}-\frac{1}{p})} \right)e^{-\mu_1 t} \| z \|_{L^{q}(\Omega)}\quad \textrm{for any }\; z\in L^q(\Omega), 
    \end{equation}
     \begin{equation}\label{SemigroupRobin2}
    \|\nabla  e^{t\Delta} z\|_{L^{p}(\Omega)} \leq C_R \left(
1+t^{-\frac{1}{2}-\frac{n}{2}(\frac{1}{q}-\frac{1}{p})} \right)e^{-\mu_1 t} \| z \|_{L^{q}(\Omega)} \quad \textrm{for any }\; z\in L^q(\Omega), 
    \end{equation}
    \begin{equation}\label{SemigroupRobin3}
    \|\nabla  e^{t\Delta} z\|_{L^{2}(\Omega)} \leq   \| \nabla z \|_{L^{2}(\Omega)} \quad \textrm{for any }\; z\in W^{1,2}(\Omega). 
    \end{equation}    
    \begin{proof}
The proofs of \eqref{SemigroupRobin1} and \eqref{SemigroupRobin2} are close to those of \cite[(1.4) and (1.5)]{WinklAggre}, respectively. Henceforth, we understand that herein it is only sufficient to show \cite[(1.9)]{WinklAggre}, i.e.,  some exponential decay estimate involving the semigroup and the principle eigenvalue for the Robin Laplacian, defined as
\begin{equation*}
\mu_1= \min_{z \in H^1(\Omega), z \neq 0} \frac{\int_\Omega |\nabla z|^2+h\int_{\partial \Omega} z^2}{\int_\Omega z^2}. 
\end{equation*}
Because \(\int_\Omega (\Delta z)z=-\int_\Omega |\nabla z|^2-h\int_{\partial \Omega} z^2\leq -\mu_1 \|z\|^2_{L^2(\Omega)}\) for all \(z \in D(\Delta_R)\), then
\[
\frac{d}{dt} \|e^{t \Delta}z\|^2_{L^2(\Omega)} = 2\int_\Omega (\Delta e^{t \Delta}z)( e^{t \Delta}z)\leq -2\mu_1 \int_\Omega (e^{t \Delta}z)^2 = -2\mu_1 \|e^{t \Delta}z\|^2_{L^2(\Omega)},
\]
for all \(z \in D(\Delta_R)\). Hence, $\|e^{t \Delta}z\|^2_{L^2(\Omega)} \leq e^{-2\mu_1 t} \|e^{0 \Delta}z\|^2_{L^2(\Omega)}$, and subsequently $ \|e^{t \Delta}z\|_{L^2(\Omega)} \leq e^{-\mu_1 t} \|z\|_{L^2(\Omega)}$; this extends to all \(z\in L^2(\Omega)\) because \(D(\Delta_R)\) is dense in $L^2(\Omega)$.

Finally, similarly to the previous derivation, we have
\[
\frac{d}{dt} \|\nabla e^{t \Delta}z\|^2_{L^2(\Omega)} = 2\int_\Omega (\Delta \nabla e^{t \Delta}z)( \nabla e^{t \Delta}z)\leq -2\mu_1 \int_\Omega |\nabla e^{t \Delta}z|^2 = -2\mu_1 \|\nabla e^{t \Delta}z\|^2_{L^2(\Omega)}\leq 0,
\]
and henceforth also \eqref{SemigroupRobin3} is established.
\end{proof}
\end{lemma}
\subsection{Parabolic and elliptic estimates}
\begin{lemma}\label{regularity}
Let $\Omega\subset\R^n$ satisfy condition in \eqref{reglocal}, $T\in(0,\infty]$, $p>n$ and $h>0$. 
\begin{itemize}
    \item [$\triangleright$] {\bf\emph{Elliptic regularity}:}  If $0\leq \psi\in C^{\delta}(\bar\Omega)$, then the solution $z\in C^{2+\delta}(\bar\Omega)$  of problem
\begin{equation*}
\begin{cases}
    -\Delta z + z=\psi  & \text{in } \Omega,\\
    z_{\nu}=-h z & \text{on } \partial\Omega,
\end{cases}
\end{equation*}
is nonnegative in $\bar{\Omega}$, and moreover it has the property that for every $q>2$ 
 there is $C_E=C_E(n,\Omega,q,h)>0$  for which $z$ is such that
\begin{equation}\label{ellipticreg}
   \lVert z \rVert_{W^{1,q}(\Omega)}\leq C_E \lVert \psi \rVert_{L^{q}(\Omega)}.
\end{equation}
Additionally, if $\psi\in L^p(\Omega)$ then
\begin{equation}\label{tau0extension} 
    z \in W^{1,\infty}(\Omega).
\end{equation}
 \item [$\triangleright$] {\bf\emph{Parabolic regularity}:} If $0\leq \psi\in C^{\delta,\frac{\delta}{2}}(\bar{\Omega}\times [0,T])$ and $z_0\in C^{2+\delta}(\Omega)$, with $z_{0\nu} = -hz_0$ on $\partial\Omega$, then the solution  $z\in C^{2+\delta,1+\frac{\delta}{2}}(\bar{\Omega}\times [0,T])$ of the problem
\begin{equation*}
\begin{cases}
    z_t= \Delta z - z  +\psi & \text{in } \Omega \times (0,T),\\
     z_\nu = -h z & \text{on } \partial\Omega \times (0,T),\\
    z(\cdot,0)=z_0 & \text{on } \Omega,
\end{cases}
\end{equation*}
is nonnegative in $\bar{\Omega}\times [0,T]$, and additionally it has the property that for every $q>1$ there exists $C_P=C_P(n,\Omega,q,h)>0$ for which $z$ is such that
\begin{equation}\label{tau1}
\int_0^t e^s \into \abs*{\Delta z(\cdot,s)}^q\,ds \leq C_{P}\quadra*{1+\int_0^t e^s \into\abs{\psi(\cdot,s)}^q\,ds} \quad \text{for all } t\in(0,T).
\end{equation}
Additionally, if  $\psi\in L^\infty((0,T);L^p(\Omega))$ then
\begin{equation}\label{tau1extension}
z \in L^\infty((0,T);W^{1,\infty}(\Omega)).
\end{equation}
\end{itemize}
\begin{proof}  
\quad $\diamond$ Existence and regularity are given in \cite[Theorem 2.4.2.6]{GrisvardBook} and \cite[Theorem 6.31]{GilbarTrudinger}, whereas nonnegativity is consequence of \cite[Theorem 2 and Hopf's Lemma, $\S$6.4]{Evans-2010-PDEs}. Let us prove \eqref{ellipticreg}.  From \cite[(3.1.9)]{LunardiBook}, it is possible to find \const{G1} (also depending on $h$) such that 
\begin{equation}\label{AAEquation}
\lVert z \rVert_{W^{2,q}(\Omega)}\leq \const{G1}\left(\lVert z \rVert_{L^q(\Omega)}+ \lVert \psi \rVert_{L^q(\Omega)}\right).
\end{equation}
On the other hand, the Gagliardo--Nirenberg inequality (\cite[page 37]{HenryBook}) entails, for $q>2$, the value $\theta=\frac{\frac{1}{2}-\frac{1}{q}}{\frac{1}{2}+\frac{2}{n}-\frac{1}{q}}\in (0,1)$ such that, by invoking the Young inequality, we have
\begin{equation*}
\lVert z \rVert_{L^{q}(\Omega)}\leq \const{G-2}\left(\lVert z \rVert_{L^{2}(\Omega)}^{1-\theta} \lVert z \rVert_{W^{2,q}(\Omega)}^\theta \right)\leq \frac{1}{2\const{G1}}\lVert z \rVert_{W^{2,q}(\Omega)}+\const{G01}\lVert z \rVert_{L^{2}(\Omega)},
\end{equation*}
which, together with bound \eqref{AAEquation}, leads to 
\begin{equation}\label{Inequ_1-G}
\lVert z \rVert_{W^{2,q}(\Omega)}\leq \const{G-a2}\left(\lVert z \rVert_{L^{2}(\Omega)}+\lVert \psi \rVert_{L^{q}(\Omega)}\right).
\end{equation}
Further, testing procedures and H\"{o}lder's inequality  provide 
\begin{equation*}
\lVert z\rVert_{L^2(\Omega)}^2\leq \int_\Omega |\nabla z|^2+h\int_{\partial \Omega} z^2 +\int_\Omega z^2=\int_\Omega z \psi \leq \lVert \psi \rVert_{L^q(\Omega)} \lVert z \rVert_{L^\frac{q}{q-1}(\Omega)}\leq \const{G_2}\lVert \psi\rVert_{L^q(\Omega)} \lVert z \rVert_{L^2(\Omega)},
\end{equation*}
i.e., 
\begin{equation}\label{Inequ_2-G}
\lVert z\rVert_{L^2(\Omega)}\leq  \const{G_2}\lVert \psi\rVert_{L^q(\Omega)},
\end{equation}
so that the claim is given by inserting relation \eqref{Inequ_2-G} into \eqref{Inequ_1-G}.
Finally, relation \eqref{tau0extension} comes from \cite[Theorem 2.4.2.7]{GrisvardBook}, which warrant that $z \in W^{2,p}(\Omega)$ and hence 
$\nabla z \in W^{1,p}(\Omega)$. Finally, since $p>n$, Sobolev embedding theorems $W^{2,p}(\Omega) \hookrightarrow W^{1,p}(\Omega) \hookrightarrow L^{\infty}(\Omega)$ imply that $z, \nabla z \in L^{\infty}(\Omega)$.

\quad $\diamond$  Existence and regularity are given in \cite[Corollary 5.1.22]{LunardiBook}, whereas nonnegativity is consequence of \cite[Theorem 9, $\S$7.1]{Evans-2010-PDEs} and  \cite[Theorem 1]{Friedman1958}.

Since from \cite[Theorem 3.1.3]{LunardiBook}, the $L^q$-realization $(A(t)z)(x)=\Delta z - z$, equipped with Robin conditions, generates an analytic semigroup on $X=L^q(\Omega)$ and $X_1=D(A)=W^{2,q}_\mathcal{\frac{\partial }{\partial \nu}}(\Omega)= \{z\in W^{2,q}(\Omega): z_\nu=-h z\,\,\textrm{on}\;\partial \Omega\}$ (see also \cite[$\S$4]{PrussScnaubelt}), we can adapt \cite[Lemma 3.6]{IshidaLankeitVigliloro-Gradient}, so to establish \eqref{tau1}.

Finally, the last claim is derived as follows. For $p>n$, we have \(- \frac{1}{2} - \frac{n}{2p} > -1 \). Subsequently, since by means of the representation formula for \( z \) we have
\begin{equation*}
z(\cdot, t) = e^{t(\Delta-1)} z_0 + \int_0^t e^{(t-s)(\Delta-1)} \psi(\cdot, s) ds \quad \text{for all } t \in (0, T),
\end{equation*}
we can invoke properties \eqref{SemigroupRobin1} and \eqref{SemigroupRobin2}, so to have for some $\const{G_x}=\const{G_x}(C_R)$
\begin{equation*}
\begin{split}
\|z(\cdot, t)\|_{W^{1,\infty}(\Omega)} & \leq e^{-t} \|e^{t\Delta} z_0\|_{W^{1,\infty}(\Omega)} + \int_0^t \|e^{(t-s)(\Delta-1)} \psi(\cdot, s)\|_{W^{1,\infty}(\Omega)} ds \\
&\leq \|e^{t\Delta} z_0\|_{W^{1,\infty}(\Omega)} + \int_0^t \|e^{(t-s)(\Delta-1)} \psi(\cdot, s)\|_{L^{\infty}(\Omega)} +\int_0^t \|\nabla e^{(t-s)(\Delta-1)} \psi(\cdot, s)\|_{L^{\infty}(\Omega)}ds\\
&\leq \const{G_x}+\const{G_x}\int_0^t \left(1+(t - s)^{- \frac{n}{2p}}\right)e^{-\mu_1(t-s)}\|\psi(\cdot, s)\|_{L^{p}(\Omega)}ds
+ \const{G_x} \int_0^t \left(1+ (t - s)^{-\frac{1}{2}- \frac{n}{2p}}\right)e^{-\mu_1(t-s)}\|\psi(\cdot, s)\|_{L^{p}(\Omega)}ds.
\end{split}
\end{equation*}
As a consequence, from $\psi\in L^\infty((0,T);L^p(\Omega))$ and the convergence of  
\[
\int_0^t \left(1+(t - s)^{- \frac{n}{2p}}\right)e^{-\mu_1(t-s)}ds \quad \textrm{and} 
\quad \int_0^t \left(1+ (t - s)^{-\frac{1}{2}- \frac{n}{2p}}\right)e^{-\mu_1(t-s)}ds,
\]
we conclude the proof. 
\end{proof}
\end{lemma}
\subsection{Trace embeddings; relates estimates}
\begin{lemma}\label{LemmaTrace}
Let $\Omega\subset\R^n$ satisfy condition in \eqref{reglocal}. Then for all $0\leq \psi\in C^{1}(\bar\Omega)$ and $\mathfrak{p}\in \{1,2\}$ there is $C_{\partial \Omega}=
C_{\partial \Omega}(n, \Omega)$ positive such that the  following holds:
\begin{equation}\label{TraceInequ}
\int_{\partial \Omega} \psi^{\mathfrak{p}}\leq C_{\partial \Omega}\left(\int_\Omega \psi^{\mathfrak{p}}  +\int_\Omega |\nabla\psi|^{\mathfrak{p}}\right).
\end{equation} 
\begin{proof}
See \cite[page 315]{BrezisBook} or \cite[Theorem 1, page 258]{Evans-2010-PDEs}.
\end{proof}
\end{lemma}
\subsection{H\"{o}lder continuity of solutions to a class of parabolic problems}
To derive the boundedness criterion presented in section $\S$\ref{SectionLocalExis}, we first establish explicit pointwise bounds for a class of parabolic problems involving positive total flux. This method, based on a Moser-type iterative argument, yields uniform estimates of the solution in some H\"{o}lder space, expressed in terms of the bounds on the initial data and the $L^1$-norm of the solution itself. (Similar results already established in the literature for zero-total flux scenarios are available in \cite{DingWinkleSmallDensity,TaoWinkParaPara}.)
\begin{lemma}\label{RecursiveInequalityLemma}
For all $k\in \N$ and  $\theta \in (0,1)$, let $\alpha_k,\beta_k \in[0,\infty)$ complying with $\alpha_k=O(k^2)$ and $\beta_k=O(k^{\frac{2\theta}{1-\theta}})$. Then there exists $L\geq1$ such that $\alpha_k\beta_k\leq L^k$, for all $k\in \N$.

Additionally, if for all \( k \in \mathbb{N}_0 \), \( M_k \in [1, \infty) \) is such that
\begin{equation*}\label{RecursiveInequality}
M_k \leq \alpha_k\beta_k M^2_{k-1} \quad \text{for all} \quad k \in \mathbb{N}, 
\end{equation*}
then
\[
M_k \leq L^{k+2^{k+1}} M^{2^k}_{0} \quad \text{for all} \quad k \in \mathbb{N}.
\]
\begin{proof}
The first assertion is evident. As to the second part, it is established by following \cite[Lemma 4.3]{WinklerExponentailDecay} (with $\theta_k=2$ and $d=0$) once $\alpha_k\beta_k\leq L^k$ is taken into account. 
\end{proof}
\end{lemma}
\begin{lemma}\label{LemmaMoserType}
Let $\Omega\subset\R^n$ satisfy condition in \eqref{reglocal}, $T\in (0,\infty]$, $a,b,c \geq 0$, $f=f(x,t) \in \left(C^1(\bar\Omega\times (0,T))\cap L^\infty((0,T);L^\infty(\Omega))\right)^n$, with $f\cdot \nu\geq 0$ for all $(x,t)\in \partial \Omega \times (0,T)$, and 
$0\leq g=g(x,t) \in  C^0(\bar\Omega\times [0,T))$. If $0 \leq \psi\in C^{2+\delta,1+\frac{\delta}{2}}(\bar\Omega \times [0,T))\cap L^\infty((0,T);L^1(\Omega))$ solves
\begin{equation*}
\begin{cases}
\psi_t= \Delta \psi +\nabla \cdot (f \psi)+a\psi-b\psi^2-c|\nabla \psi|^2 & \textrm{in }\; \Omega \times (0,T)\\
\psi_\nu +g\psi= 0 & \textrm{on }\; \partial \Omega \times (0,T),\\
\end{cases}
\end{equation*}
and satisfies $\psi_\nu(x,0)+g(x,0)\psi(x,0)=0$ for all $x\in \partial \Omega$, then it holds that $\psi \in L^\infty((0,T);C^{1+\delta}(\bar{\Omega}))$.
\begin{proof}
For any $k\in \N_0$, let us set $p_k=2^k\in [1,\infty)$. Standard testing procedures 
allow to write for all $t\in (0,T)$ and for $k \geq 1$
\begin{equation}\label{Estimate1G}
\begin{split}
&\frac{d}{dt} \int_\Omega \psi^{p_k} \leq  -p_k(p_k-1) \int_\Omega \psi^{{p_k}-2} |\nabla \psi|^2 + {p_k}\int_{\partial \Omega} (f\cdot \nu)\psi^{p_k} 
- {p_k} ({p_k}-1) \int_\Omega \psi^{{p_k}-1} \nabla \psi \cdot f +a {p_k} \int_\Omega \psi^{p_k}\\
&\leq  -{p_k}({p_k}-1) \int_\Omega \psi^{{p_k}-2} |\nabla \psi|^2 + \const{xyz}{p_k} \int_{\partial \Omega} \psi^{p_k}  +\frac{{p_k}({p_k}-1)}{4}\int_\Omega \psi^{{p_k}-2}|\nabla \psi|^2+\const{xyz}^2{p_k}({p_k}-1)\int_\Omega \psi^{p_k}+a {p_k} \int_\Omega \psi^{p_k}\\
& \leq -\frac{3{p_k}({p_k}-1)}{4}\int_\Omega \psi^{{p_k}-2} |\nabla \psi|^2+{p_k}\left(C_{\partial \Omega}\const{xyz}+a+\const{xyz}^2({p_k}-1)\right)\int_\Omega \psi^{{p_k}}+{p_k}^2 \const{xyz}C_{\partial \Omega} \int_\Omega \psi^{{p_k}-1} |\nabla \psi|,
 \end{split}
\end{equation}
where we have used Young's inequality, $||f(\cdot,t)||_{L^\infty(\Omega)}\leq \const{xyz}$ for all $t\in (0,T)$ and introduced the trace constant of relation \eqref{TraceInequ} (with $\mathfrak{p}=1$). Successively we make use again of the Young inequality to rearrange the term involving $\int_\Omega \psi^{{p_k}-1} |\nabla \psi|$; relation \eqref{Estimate1G} becomes 
\begin{equation}\label{Estimate2G}
\begin{split}
\frac{d}{dt} \int_\Omega \psi^{p_k} 
 \leq -2\frac{{p_k}-1}{p_k}\int_\Omega |\nabla \psi^\frac{p_k}{2}|^2+{p_k}\left(C_{\partial \Omega}\const{xyz}+a+\const{xyz}^2({p_k}-1)+\frac{{p_k}^2\const{xyz}^2 C_{\partial \Omega}^2}{p_k-1}\right)\int_\Omega \psi^{{p_k}}\quad \textrm{on }\, (0,T).
 \end{split}
\end{equation}
Now we introduce $r \in \left(2,\frac{2n}{(n-2)_+}\right)$ and $\theta=\frac{1-\frac{1}{r}}{\frac{1}{2}+\frac{1}{n}}\in (0,1)$ so to be able to write thanks to the Gagliardo--Nirenberg inequality
\begin{equation}\label{GNp_k}
||\psi^\frac{p_k}{2}||_{L^r(\Omega)}^2\leq \const{SGKY} \left(||\nabla \psi^\frac{p_k}{2}||_{L^2(\Omega)}^{2\theta}||\psi^\frac{p_k}{2}||_{L^1(\Omega)}^{2(1-\theta)}+||\psi^\frac{p_k}{2}||_{L^1(\Omega)}^2\right)\quad \textrm{for all }\,t<T,
\end{equation}
for some $k$-independent $\const{SGKY}$. 
On the other hand, H\"{o}lder's inequality provides
\begin{equation}\label{Holder}
\int_\Omega \psi^{p_k}\leq |\Omega|^\frac{r-2}{r}||\psi^\frac{p_k}{2}||_{L^r(\Omega)}^2\quad \textrm{for all }\,t<T,
\end{equation}
so that for 
\begin{equation*}\label{HkDefinition}
H_k:=\alpha_k\beta_k,\quad \textrm{with }\,
\begin{cases}
\alpha_k=p_k\const{xyz}C_{\partial \Omega}+ap_k+{\const{xyz}}^2p_k(p_k-1)+\frac{p_k^3 \const{xyz}^2C_{\partial \Omega}^2}{p_k-1}+1,\\
\beta_k=\const{PDDP}  \left(\left(\frac{p_k-1}{p_k\alpha_k}\right)^\frac{\theta}{\theta-1}+1\right),
\end{cases}
\end{equation*}
relation \eqref{GNp_k} supported by the Young inequality implies that estimate \eqref{Estimate2G}, also applying \eqref{Holder}, is for some $k$-independent constant $\const{PDDP}$ reduced to
\begin{equation*}
\frac{d}{d t}\int_\Omega \psi^{p_k}+\int_\Omega \psi^{p_k}\leq  H_k \left(\int_\Omega \psi^\frac{p_k}{2}\right)^2\quad \textrm{for all }\,t\in(0,T).
\end{equation*}  
Successively, since $\alpha_k$ and $\beta_k$ fulfill assumptions in Lemma \ref{RecursiveInequalityLemma}, the sequence
\begin{equation}\label{DefinitionSequelMk0}
M_k(T):=\sup\bigg\{1,\sup_{t\in (0,T)}\int_\Omega \psi^{p_k}\bigg\}\quad \textrm{for all } k\in \N_0,
\end{equation}
obeys for some $L>1$ this initial value problem 
\begin{equation}\label{ODI_M_k}
\frac{d}{d t} \int_\Omega \psi^{p_k} \leq -\int_\Omega \psi^{p_k} +L^kM^2_{k-1}(T) \quad \textrm{ on } (0,T)\;\;\textrm{and }\;\;
\int_\Omega \psi(x,0)^{p_k}dx=:\int_\Omega \psi^{p_k}_0(x)dx.
\end{equation}
Since from hypothesis $M_0(T)=\sup\left\{1,\sup_{t \in (0,T)}\int_\Omega \psi(x,t)dx\right\}$ is
finite and time independent, by integrating the above ODI we have that 
\begin{equation*}
M_1(T)\leq \max\Big\{\int_\Omega \psi^2_0(x)dx,LM_{0}^{2}\Big\},
\end{equation*}
i.e., $M_1(T)$ is uniformly bounded on $(0,T)$ as well. Definitely, by induction, we deduce that
actually relation \eqref{DefinitionSequelMk0} becomes
\begin{equation}\label{DefinitionSequelMk}
M_k:=\max\bigg\{1,\sup_{t\in (0,T)}\int_\Omega \psi^{p_k}\bigg\}\quad \textrm{for all } k\in \N_0,
\end{equation}
where $M_k$ fulfills thanks to \eqref{ODI_M_k}
\begin{equation*}
M_k\leq \max\Big\{\int_\Omega \psi_0^{p_k}(x)dx,L^kM_{k-1}^{2}\Big\}  \quad \textrm{ for all } k\in \N.
\end{equation*}
In these circumstances,  we have to distinguish two cases. 
\begin{itemize}
	\item  [$\triangleright$]If there exists a subsequence of natural numbers $k_j$ such that $k_j\nearrow \infty$ as $j \nearrow \infty$ and 
	$$M_{k_j}\leq  \int_\Omega \psi^{p_{k_j}}_0(x)dx\quad \textrm{for all } j\in \N,$$	
then by virtue of \eqref{DefinitionSequelMk} we have
\begin{equation}\label{BoundIterative1}
\begin{split}
\lVert \psi(\cdot, t)\rVert_{L^\infty(\Omega)}&:= \lim_{j\rightarrow \infty}\bigg(\int_\Omega \psi^{p_{{k_j}}}\bigg)^\frac{1}{p_{{k_j}}} 
\leq \limsup_{j\rightarrow \infty} M_{{k_j}}^{\frac{1}{p_{{k_j}}}}
 \leq \limsup_{j\rightarrow \infty} \bigg(\int_\Omega \psi_0^{p_{{k_j}}}(x)dx\bigg)^\frac{1}{p_{{k_j}}}\leq \lVert \psi_0 \rVert_{L^\infty(\Omega)} \quad \textrm{for all } t\in  (0,T).
\end{split}
\end{equation}
	\item  [$\triangleright$] Conversely, if such a sequence does not exist, we hence arrive at 
	\begin{equation*}
		M_k\leq L^{k} M_{k-1}^{2}\quad \textrm{ for all } k\in \N,
	\end{equation*}
and an application of Lemma \ref{RecursiveInequalityLemma} infers
\begin{equation*}
	M_k\leq L^{k+2^{k+1}} M_{0}^{ 2^k}  \quad \textrm{ for all } k\in \N,
\end{equation*}
and in turn 
\begin{equation}\label{BoundIterative2}
\begin{split}
\lVert \psi (\cdot, t)\rVert_{L^\infty(\Omega)}&:=
\lim_{k\rightarrow \infty} \bigg(\int_\Omega \psi^{p_{k}}\bigg)^\frac{1}{p_{k}} 
\leq \limsup_{k\rightarrow \infty} M_{k}^{\frac{1}{p_{k}}}
\leq \limsup_{k\rightarrow \infty} L^\frac{k+2^{k+1}}{p_k} M_{0}^{{\frac{2^{k}}{p_k}}}=L^2 M_0\quad \textrm{on }  (0,T).
\end{split}
\end{equation}
\end{itemize}
From the last estimates \eqref{BoundIterative1} and \eqref{BoundIterative2} we get
$ 
\lVert \psi (\cdot, t)\rVert_{L^\infty(\Omega)}\leq \max\{\lVert \psi_0 \rVert_{L^\infty(\Omega)},L^2 M_0\}
$, i.e., we have 
$\psi\in L^\infty((0,T);L^\infty(\Omega)).$
Finally, we have boundedness of $\lVert \psi(\cdot,t) \rVert_{C^{1+\delta}(\bar{\Omega})}$ on $[0,T]$, once \cite[Theorem ~1.2]{lieberman_paper} is applied.
\end{proof}
\end{lemma}
\section{Local existence, dichotomy and boundedness criterion}\label{SectionLocalExis}
Let us establish existence and uniqueness of local classical solutions to problem \eqref{problem}. We provide the proof because, although the reasoning for the two scenarios is similar, the literature only addresses the case with Neumann boundary conditions and without gradient nonlinearities. 
\begin{lemma}\label{theoremExistence}
Let the hypotheses in \eqref{reglocal} be fulfilled.  Additionally, let $a,b, c \geq 0$,
$\chi,h >0$, $\tau \in \{0,1\}$ and $\alpha \in[0,1]$. Then for every 
 \begin{equation*}
 u_0, \tau v_0: \bar{\Omega}  \rightarrow \R^+ , \textrm{ with } u_0,  \tau v_0 \in C^{2+\delta}(\bar\Omega) \quad \textrm{complying with} \quad u_{0\nu}=(\alpha-1)\chi h u_0v_0 \quad \textrm{ and } \quad \tau v_{0\nu}=-h \tau v_0 \quad  \textrm{ on }\partial \Omega
\end{equation*}
there exist $\TM\in (0,\infty]$ and a unique couple of functions $(u,v$), with 
\begin{equation*}
(u,v)\in C^{2+\delta,1+\frac{\delta}{2}}( \Bar{\Omega} \times [0, \TM))\times C^{2+\delta,\tau+\frac{\delta}{2}}( \Bar{\Omega} \times [0, \TM)),
\end{equation*}
solving problem \eqref{problem} and nonnegative in $\bar{\Omega}\times [0,\TM)$. Additionally,
 \begin{equation}\label{dictomyCriteC2+del} 
 \text{if} \quad \TM<\infty \quad \text{then} \quad \limsup_{t \to \TM} \left(\|u(\cdot,t)\|_{C^{2+\delta}(\bar\Omega)}+\|v(\cdot,t)\|_{C^{2+\delta}(\bar\Omega)}\right)=\infty.
 \end{equation}
\begin{proof}
Let us start with the \textit{existence} issue.  For any $R>0$,
let us consider for some $0<T\leq 1$, to be precised later on,  the closed convex
subset
$$S_{T}=\{0\le u \in C^{1,\frac{\delta}{2}}(\bar{\Omega}\times [0,T]):  \|u(\cdot,t)-u_0\|_{L^\infty(\Omega)}\le R,\, \; \text{for\ all}\
t\in[0,T]\}.$$
Once an element $\tilde{u}$ of $S_T$ is picked, Lemma \ref{regularity} provides the solution $v$ to problem 
 \begin{equation*}\label{2.2}
		\begin{cases}
	\tau v_t-	\Delta v+ v=\tilde{u}&\text{in}\ \Omega\times(0,T),\\
		v_\nu=-hv&\text{on}\ \partial\Omega\times(0,T),\\
		\tau v(x,0)=\tau v_0(x),
	\end{cases}
\end{equation*}
which is, thanks to $\tau v_0\in C_\nu^{2+\delta}(\bar{\Omega})$ and  elliptic and parabolic regularity results for oblique derivative problems (\cite[Theorem 6.30]{GilbarTrudinger} for $\tau=0$, and \cite[Corollary 5.1.22]{LunardiBook}, \cite[Theorem IV.5.3]{LSUBookInequality} for $\tau=1$), such that 
\begin{equation*}
v\in C^{2+\delta, \tau+\frac{\delta}{2}}(\bar\Omega \times [0,T])
\quad \textrm{and precisely}\quad \sup_{t\in [0,T]} \lVert v(\cdot, t)\rVert_{C^{2+\delta}(\bar\Omega)}\leq H,
\end{equation*}
where $H=H(\tilde{u},h,\tau v_0)>0$ is estimated by $\tilde{u}$, so essentially, $H=H(R)$.
Now, taken into account the properties of $v$ and $\nabla v$ on $\bar{\Omega}\times [0,T]$,  the problem
\begin{equation}\label{2.1}
	\begin{cases}
		u_t=\Delta u -\chi  \nabla u  \cdot \nabla v-\chi u \Delta v +  a u-b u^2 - c\abs{\nabla u}^2
		&\text{in}\ \Omega\times(0,T),\\
		u_\nu=(\alpha-1)h \chi u v&\text{on}\ \partial\Omega\times(0,T),\\
		u(x,0)=u_0(x)&x\in\overline{\Omega},
	\end{cases}
\end{equation}
is, with an evident choice of $\mathcal{A}$ and $\varphi$, a general system of the form
\begin{equation*}
	\begin{cases}
		u_t=\mathcal{A}u +\varphi(x,t,u,\nabla u)
		&\text{in}\ \Omega\times(0,T),\\
	\mathcal{B}_1u=u_\nu-(\alpha-1)h \chi u v=0&\text{on}\ \partial\Omega\times(0,T),\\
		u(x,0)=u_0(x)&x\in\overline{\Omega}.
	\end{cases}
\end{equation*}
Since the above problem complies with \cite[(7.3.5)]{LunardiBook},  by recalling that $u_0\in C^{2+\delta}_\nu(\bar\Omega)$,  \cite[Proposition 7.3.3]{LunardiBook} and \cite[Theorem V 6.1]{LSUBookInequality} (see also \cite[Theorem 8.5.4]{LunardiBook}) ensure the existence of some $0<T< 1$, as above, such that problem \eqref{2.1} has a unique solution 
\begin{equation*}
u\in C^{2+\delta,1+\frac{\delta}{2}}(\bar\Omega\times [0,T]).
\end{equation*} 
In particular,  this produces some positive constant $K=K(\tilde{u},v,\nabla v,h,\alpha,u_0)$, and henceforth $K=K(R)$, with the property that
$$|u(x,t)-u_0(x)|\le Kt^{1+\tfrac{\delta}{2}} \quad \textrm{for all} \quad(x,t)\in\Omega\times(0,T), \quad \textrm{or}\quad \max_{t\in[0,T]}\,\|u(\cdot,t)-u_0\|_{L^\infty(\Omega)}\le KT^{1+\frac{\delta}{2}}. $$
In this way, we deduce 
that 
$$ \textrm{for} \quad T\leq \left(\frac{R}{K}\right)^\frac{2}{2+\delta}, \quad\|u(\cdot,t)-u_0\|_{L^\infty(\Omega)}\le R\;\quad
\textrm{for all} \quad t\in[0,T]. $$
Moreover, since $\underline{u}\equiv 0$ is a subsolution of the first equation in \eqref{2.1}, the parabolic comparison principle (see \cite[Theorem 8.2.]{CrandallEtAl_User1992} supported by \cite[Theorem 1]{Friedman1958}) warrants the nonnegativity of $u$ on $\Omega \times (0,T)$. So, for values of $T$ as above, the map  $\Phi (\tilde{u})=u$ where $u$ solves problem \eqref{2.1}, is such that $\Phi(S_T) \subset S_T$ and $\Phi$ is compact, because the
 \cite[Ascoli--Arzel\`a Theorem 4.25]{BrezisBook} implies that the natural embedding
of $C^{2+\delta,\tau+\frac{\delta}{2}}(\bar{\Omega}\times [0,T]))$ into $C^{1,\frac{\delta}{2}}(\bar{\Omega}\times [0,T]))$ is a compact linear operator, and henceforth the Schauder's fixed point theorem asserts the existence of a fixed point $u$ for $\Phi$. 

As to \textit{uniqueness}, let $(u_1,v_1)$ and $(u_2,v_2)$ be two different nonnegative classical solutions of problem \eqref{problem} in $\Omega\times (0,T)$ with the same initial data $u_1(\cdot,0)=u_2(\cdot,0)=u_0(x)$ and $\tau v_1(\cdot,0)=\tau v_2(\cdot,0)=\tau v_0(x)$. By confining our attention to the case $\tau=0$ (being the case $\tau=1$ similar, as we will indicate later), we thus have for $i=1,2$  these problems:
\begin{equation}\label{3.2_Bis}
	\begin{cases}
		-\Delta v_i+ v_i=u_i&\text{in}\ \Omega\times(0,T),\\
		(v_i)_\nu=-hv_i&\text{on}\ \partial\Omega\times(0,T),
	\end{cases}
\end{equation}
and
\begin{equation}\label{3.1_BisA}
	\begin{cases}
  (u_i)_t=\Delta u_i
		-\chi \nabla \cdot (u_i \nabla v_i)
		+a u_i-b u_i^2-c|\nabla u_i|^2 
 & \text{in}\ \Omega\times(0,T), 
\\
		(u_i)_\nu=(\alpha-1)\chi h u_i v_i&\text{on}\ \partial\Omega\times(0,T),
  \\
		u_i(x,0)=u_0(x)&x\in\overline{\Omega}.
	\end{cases}
\end{equation}
From the $C^2$-regularity of $u_1, u_2, v_1$and $v_2$, we can set 
\begin{equation}\label{ConstantsC1-2-3} 
\begin{split}
	C_1=
		 & \max\{\|u_1\|_{L^\infty(\Omega\times(0,T))},
	\|u_2\|_{L^\infty(\Omega\times(0,T))}, \lVert v_1\rVert_{L^\infty(\Omega \times (0,T))}, \lVert v_2\rVert_{L^\infty(\Omega \times (0,T))}, \\ & \quad \quad  \lVert \nabla u_1\rVert_{L^\infty(\Omega \times (0,T))},\lVert \nabla u_2\rVert_{L^\infty(\Omega \times (0,T))}, \lVert \nabla v_1\rVert_{L^\infty(\Omega \times (0,T))}\}.
	\end{split}
\end{equation}
Thereafter, the Mean Value Theorem implies  
\begin{equation}\label{MeanValueTheroem}  
|u_1^2-u_2^2| \leq 2 C_1|u_1-u_2| \quad \textrm{ and } \quad ||\nabla u_1|^2-|\nabla u_2|^2|\leq 2 C_1||\nabla u_1|-|\nabla u_2|| \quad \text{ in}\ \Omega\times(0,T).
\end{equation} 
%
Under such circumstances, by considering problems \eqref{3.2_Bis}, we have that $V:=v_1-v_2$ solves
\begin{equation}\label{ProblElliptic}
\begin{cases}
-\Delta V+V=U:=u_1-u_2&\text{in}\ \Omega\times(0,T),\\
V_\nu=-h V & \text{on}\ \partial\Omega\times(0,T).
\end{cases}
\end{equation}
 In light of this, through Young's inequality, testing against $V$ the mentioned equation provides the following estimate
\begin{align*}
h\int_{\partial \Omega} V^2	+\int_\Omega|\nabla V|^2&+\int_\Omega V^2=\int_\Omega UV 
	\le\frac{1}{2}\int_\Omega V^2
	+\const{Giu-2}\int_\Omega U^2 \quad \textrm{on }\,\,(0,T),
\end{align*}
inferring 
\begin{equation}\label{3.4}
	\int_\Omega|\nabla V|^2+\int_\Omega V^2\le
	\const{Giu-2Bis}\int_\Omega U^2\quad \textrm{for all }\, t\in(0,T).
\end{equation}
On the other hand, from \eqref{3.1_BisA} some computations show that $U=u_1-u_2$ solves
\begin{equation}\label{3.1_Bis}
	\begin{cases}
  U_t=\Delta U
		-\chi \nabla \cdot (U \nabla v_1+u_2\nabla V)
		+a U -b (u_1^2-u_2^2)-c(|\nabla u_1|^2-|\nabla u_2|^2) 
 & \text{in}\ \Omega\times(0,T), 
\\
		U_\nu=(\alpha-1)\chi h\left( u_1 v_1-u_2v_2\right)&\text{on}\ \partial\Omega\times(0,T),
  \\
		U(x,0)=0&x\in\overline{\Omega}.
	\end{cases}
\end{equation}
In this way, from problem \eqref{3.1_Bis} we get
\begin{equation*}
\begin{split}
\frac{1}{2}\frac{d}{dt}\int_\Omega U^2+\int_\Omega |\nabla U|^2&=(\alpha-1)\chi h \int_{\partial \Omega}U(u_1v_1-u_2v_2)-\chi \int_{\partial \Omega}U^2 v_{1\nu}-\chi \int_{\partial \Omega}U u_2 V_{\nu}\\
	& \quad +\chi \int_\Omega U \nabla U \cdot \nabla v_1+\chi \int_\Omega u_2 \nabla U \cdot \nabla V +\int_\Omega a U^2\\
	& \quad -b\int_\Omega  U (u_1^2-u_2^2) -c\int_\Omega U (|\nabla u_1|^2-|\nabla u_2|^2) \quad \textrm{ for all } \,t\in (0,T),
	\end{split}
\end{equation*}
or also, taken into account \eqref{3.2_Bis}, \eqref{3.1_BisA},  \eqref{ConstantsC1-2-3},   \eqref{MeanValueTheroem}, \eqref{ProblElliptic} and the inequality $|\nabla u_1|-|\nabla u_2|\leq |\nabla (u_1-u_2)|=|\nabla U|$, 
\begin{equation}\label{MainForUniqueness}  
\begin{split}
\frac{1}{2}\frac{d}{dt}\int_\Omega U^2+\int_\Omega |\nabla U|^2&\leq C_2\int_{\partial \Omega}U^2+C_3 \int_{\partial \Omega}|U  V| +C_4\int_\Omega |U \nabla U| +C_5\int_\Omega |\nabla U \cdot \nabla V| +C_6\int_\Omega  U^2,
	\end{split}
\end{equation}
as well valid for all $t\in (0,T)$.

In turn we invoke inequality \eqref{TraceInequ} of Lemma \ref{LemmaTrace} (with $\mathfrak{p}=2$), supported by Young's inequality,  to control the integrals at the right hand side of \eqref{MainForUniqueness}; we have for $C_7$--$C_{11}$ also possibly depending on $C_{\partial \Omega}$,  
\begin{equation}\label{IneqA-A}
\begin{cases}
C_2\displaystyle \int_{\partial \Omega}U^2+C_3 \int_{\partial \Omega}|U  V|\leq C_7\int_\Omega U^2+\frac{1}{3}\int_\Omega |\nabla U|^2+C_8\int_\Omega V^2+C_9\int_\Omega |\nabla V|^2\\
C_4\displaystyle\int_\Omega |U \nabla U| +C_5\int_\Omega |\nabla U \cdot \nabla V|\leq \frac{1}{3}\int_\Omega |\nabla U|^2+C_{10}\int_\Omega U^2+\frac{1}{3}\int_\Omega |\nabla U|^2+C_{11}\int_\Omega |\nabla V|^2.
\end{cases}
\end{equation}
Finally, by plugging  \eqref{IneqA-A} into \eqref{MainForUniqueness} yields thanks to relation \eqref{3.4}   
a suitable positive constant $\tilde{C}=\tilde{C}(T,C_{\partial \Omega})$ such that 
$$\frac{d}{dt}\int_\Omega U^2  \le\tilde{C}\int_\Omega U^2,\qquad t\in(0,T).$$ From the initial
condition $\int_\Omega U^2(x,0)\,dx=0$ (recall $U(x,0)=0$ in $\bar{\Omega}$) we establish $U=u_1-u_2=0$ on $\Omega \times (0,T)$ and  consequently, from problem \eqref{3.2_Bis}, we manifestly get $v_1-v_2=0$, as well on $\Omega \times (0,T)$; the proof is concluded. (For $\tau=1$, the proof is similar but the reasoning is based on the evolutive behavior of $\frac{1}{2}\int_\Omega \left(U^2+V^2\right).$)

Taking $T$ as initial time and $u(\cdot, T)$ as the initial condition, the above reasoning would provide a solution $(\hat{u},\hat{v})$ defined on $\bar{\Omega}\times [T,\hat{T}]$, for some $\hat{T}>0$, which by uniqueness would be the prolongation of $(u,v)$ (exactly, from $\bar{\Omega}\times [0,T]$ to  $\bar{\Omega}\times [0,\hat{T}])$. This procedure can be repeated up to construct a maximal interval time $[0,\TM)$ of existence, in the sense that either $\TM=\infty$, or if $\TM<\infty$ no solution belonging to $C^{2+\delta,\tau+\frac{\delta}{2}}(\bar{\Omega}\times [0,\TM])$ may exist and henceforth relation \eqref{dictomyCriteC2+del} has to be fulfilled. 
\end{proof}
\end{lemma}
From now on with $(u,v)$ we will refer to the unique local solution to model \eqref{problem} defined on $\Omega \times (0,\TM)$ and  provided by Lemma \ref{theoremExistence}. 
As seen in such a result, relation \eqref{dictomyCriteC2+del} holds and in this case the solution is said to blow-up at $\TM$ in $C^{2+\delta}(\Omega)$-norm; due to the quadratic growth in the gradient term of the logistic, such a blow-up implies also the explosion in $L^\infty(\Omega)$-norm.
 We precisely have this 
\begin{lemma}[Extension and boundedness criteria] \label{ExtensionLemma} If $u\in L^\infty((0,\TM);C^{1+\delta}(\bar\Omega))$, then we have $\TM=\infty$ and in particular $u\in L^\infty((0,\infty);L^\infty(\Omega))$.
\begin{proof}
By contradiction let $\TM$ be finite. Naturally, it is sufficient to show that the uniform-in-time boundedness of $u$ in $C^{1+\delta}(\bar\Omega)$ entails finiteness of $\lVert u(\cdot,t) \rVert_{C^{2+\delta}(\bar{\Omega})}+\lVert v(\cdot,t) \rVert_{C^{2+\delta}(\bar{\Omega})}$ on $[0,\TM].$ In fact, from \eqref{dictomyCriteC2+del}, we would have an inconsistency. But, from this regularity of $u$ and the relative bound, bootstrap  arguments already developed in Lemma \ref{theoremExistence}, in conjunction with the regularity of the initial data $(u_0,\tau v_0)$ for model \eqref{problem} give 
$$
\sup_{t\in [0,\TM]}\left(\lVert u(\cdot,t) \rVert_{C^{2+\delta}(\bar{\Omega})}+\lVert v(\cdot,t) \rVert_{C^{2+\delta}(\bar{\Omega})}\right)<\infty.
$$
\end{proof}
\end{lemma} 
\section{A priori estimates and proof of Theorem \ref{theoremlocal}}
Since the uniform-in-time boundedness of $u$ is implied whenever $u\in L^\infty((0,\TM);L^p(\Omega))$ for some $p>1$, here under we dedicate to the derivation of some \textit{a priori} integral estimates. The first step toward such a conclusion is deriving boundedness of $\int_\Omega u$ for all $t \in (0,\TM)$. 
\subsection{Uniform-in-time boundedness of the mass for $\tau=0$} 
\begin{lemma}\label{BoundednessMass}
Let $\tau=0$, $\chi, a, h >0$, $\alpha \in [0,1]$, $C_{\partial \Omega}$ be the trace constant in Lemma \ref{LemmaTrace} and either 
\begin{equation}\label{Estimate_bLemma}
\begin{cases}
b> c+ \frac{(\chi \alpha)^2 h C_{\partial \Omega}}{16 c}  & \text{ if } 0<c<\frac{\chi \alpha}{4} \sqrt{h C_{\partial \Omega}},\\
b>\frac{\chi \alpha}{2} \sqrt{h C_{\partial \Omega}} & \text{ if } c\geq\frac{\chi \alpha}{4} \sqrt{h C_{\partial \Omega}},
\end{cases}
\end{equation}
or 
\begin{equation}\label{BagTao0Lem}
\alpha \chi  < \min\{b, 4c\}. 
\end{equation}
Then $u\in L^\infty((0,\TM);L^1(\Omega))$, in the sense that there exists $m>0$ such that 
\begin{equation*}
\int_\Omega u \leq m \quad \textrm{for all } t \in (0, \TM).
\end{equation*}
\begin{proof}
Let us proceed to integrate the first equation of problem \eqref{problem} on $\Omega$. Through the application of 
the divergence theorem, in conjunction with the Robin boundary conditions and the invocation of the Young inequality (with $\epsilon_1>0$ chosen later), we arrive at the following estimate
\begin{equation}\label{estimate}
\begin{split}
\frac{d}{dt} \int_\Omega u &= \chi \alpha h \int_{\partial \Omega} uv + a \int_\Omega u - b \int_\Omega u^2 - c \int_\Omega  |\nabla u|^2\\
& \leq \frac{\chi^2 \alpha^2 h}{4 \epsilon_1} \int_{\partial \Omega} u^2 + \epsilon_1 h \int_{\partial \Omega} v^2 + a \int_\Omega u - b \int_\Omega u^2 - c \int_\Omega  |\nabla u|^2 \quad \textrm{on } (0, \TM).
\end{split}
\end{equation}
On the other hand, by multiplying the second equation of \eqref{problem} for $\epsilon_1 v$, integrating on $\Omega$ and applying Young's inequality (with $\epsilon_2>0$ chosen later), we obtain
\begin{equation*}
\begin{split}
&\epsilon_1  \int_\Omega |\nabla v|^2 + \epsilon_1 h \int_{\partial \Omega} v^2 + \epsilon_1   \int_\Omega v^2
= \epsilon_1  \int_\Omega uv \leq \epsilon_2
\int_\Omega u^2 + \frac{\epsilon_1^2}{4 \epsilon_2}  \int_\Omega v^2 \quad \textrm{for all } t \in (0, \TM),
\end{split}
\end{equation*}
which implies
\begin{equation*}
\epsilon_1 h  \int_{\partial \Omega} v^2 \leq \epsilon_2 \int_\Omega u^2 - \epsilon_1 \left(1-\frac{\epsilon_1}{4 \epsilon_2}\right) \int_\Omega v^2 \quad \textrm{on } (0, \TM).
\end{equation*}
By inserting the previous estimate in \eqref{estimate}, applying the trace inequality to $\frac{\chi^2 \alpha^2 h}{4 \epsilon_1} \int_{\partial \Omega} u^2$ given in Lemma \ref{LemmaTrace} (with $\psi=u$, $\mathfrak{p}=2$), this allows us to deduce
for all $t \in (0, \TM)$
\begin{equation*}\label{estimateM}
\begin{split}
\frac{d}{dt} \int_\Omega u &\leq  a \int_\Omega u + \left(\frac{\chi^2 \alpha^2 h C_{\partial \Omega}}{4 \epsilon_1} + \epsilon_2 - b \right) \int_\Omega u^2 
+ \left(\frac{\chi^2 \alpha^2 h C_{\partial \Omega}}{4 \epsilon_1}- c \right) \int_\Omega  |\nabla u|^2 - \epsilon_1 \left(1-\frac{\epsilon_1}{4 \epsilon_2}\right) \int_\Omega v^2.  
\end{split}
\end{equation*}
If $\epsilon_1$ and $\epsilon_2$ are picked with the properties that 
\begin{equation}\label{YuyaConditionDelEp1Ell}
\frac{\chi^2 \alpha^2 h C_{\partial \Omega}}{4 \epsilon_1} + \epsilon_2 - b < 0, \quad \frac{\chi^2 \alpha^2 h C_{\partial \Omega}}{4 \epsilon_1} - c \leq 0 \quad \text{and} \quad 1 - \frac{\epsilon_1}{4 \epsilon_2} > 0, 
\end{equation}
then we can apply Young's inequality to obtain constants \const{yy1} and \const{yy2} such that
\[
\frac{d}{dt} \int_{\Omega} u \leq -\const{yy1} \int_{\Omega} u + \const{yy2} \quad \text{on} \quad (0, \TM),
\]
which, together with $\int_{\Omega} u(x,0)\,dx=\int_{\Omega} u_0(x)\,dx$, implies the claim.
It remains to show the existence of $\epsilon_1$ and $\epsilon_2$ such that 
\eqref {YuyaConditionDelEp1Ell} hold. To this end, we first observe that by assumptions \eqref{Estimate_bLemma} one has 
\begin{equation*}
\frac{\alpha^2 \chi^2 h C_{\partial \Omega}}{4c} < 2b + \sqrt{4b^2 - \alpha^2 \chi^2 h C_{\partial \Omega}}. 
\end{equation*}
Indeed, if \( c < \frac{\alpha \chi \sqrt{h C_{\partial \Omega}}}{4} \), then
\[
2b + \sqrt{4b^2 - \alpha^2 \chi^2 h C_{\partial \Omega}} > 2c + \frac{\alpha^2 \chi^2 h C_{\partial \Omega}}{8c} + \sqrt{\left( \frac{\alpha^2 \chi^2 h C_{\partial \Omega}}{8c} - 2c \right)^2} = \frac{\alpha^2 \chi^2 h C_{\partial \Omega}}{4c},
\]
whereas if $c \geq \frac{\alpha \chi \sqrt{h C_{\partial \Omega}}}{4}$, it is seen that 
\begin{equation*}\label{YuyaIneqFinal_1}
2b + \sqrt{4b^2 - \alpha^2 \chi^2 h C_{\partial \Omega}} > \alpha \chi \sqrt{h C_{\partial \Omega}} \geq \frac{\alpha^2 \chi^2 h C_{\partial \Omega}}{4c}.
\end{equation*}
Thus, we can take $\epsilon_1 > 0$ such that
\begin{equation}\label{YuyaIneqFinal_0}
2b - \sqrt{4b^2 - \alpha^2 \chi^2 h C_{\partial \Omega}} < \epsilon_1 < 2b + \sqrt{4b^2 - \alpha^2 \chi^2 h C_{\partial \Omega}} \quad \text{and} \quad \epsilon_1 \geq \frac{\alpha^2 \chi^2 h C_{\partial \Omega}}{4c}. 
\end{equation}
By virtue of the first relation in \eqref{YuyaIneqFinal_0}, we have \( \epsilon_1^2 - 4b\epsilon_1 + \alpha^2 \chi^2 h C_{\partial \Omega} < 0 \), which can be rewritten as
\[
\frac{\epsilon_1}{4} < b - \frac{\alpha^2 \chi^2 h C_{\partial \Omega}}{4 \epsilon_1}.
\]
Therefore, we can choose \( \epsilon_2 > 0 \) such that
\[
\frac{\epsilon_1}{4} < \epsilon_2 < b - \frac{\alpha^2 \chi^2 h C_{\partial \Omega}}{4 \epsilon_1}.
\]
The last inequality, along with the second relation in \eqref{YuyaIneqFinal_0}, guarantees \eqref{YuyaConditionDelEp1Ell}, which leads to the conclusion.

On the other hand, we can achieve the aim reasoning in a different way, and precisely 
making use of \eqref{BagTao0Lem}. More precisely,  by adding the term $\int_\Omega u$ to both sides of equality \eqref{estimate}, we get 
\begin{equation}\label{Bag00}
\frac{d}{dt} \int_\Omega u  + \int_\Omega u  = -c \int_\Omega |\nabla u|^2 - b \int_\Omega u^2  + (a + 1) \int_\Omega u  + \chi \alpha h \int_{\partial \Omega} uv \quad \textrm{on } (0,\TM).
\end{equation}
Successively, for $\alpha \chi \leq \beta_1 < \min \{b, 4c\}$, we multiply both sides of the second equation in problem \eqref{problem} by $\beta_1(u + v)$ and integrate over $\Omega$; we obtain on $(0,\TM)$ 
\begin{equation}\label{Bag01}
0 = \int_\Omega \beta_1(u + v) \left( \Delta v + u - v \right) = -\beta_1 \int_\Omega (\nabla u \cdot \nabla v)  - \beta_1 \int_\Omega |\nabla v|^2 - \beta_1 h \int_{\partial \Omega} (uv + v^2) 
+ \beta_1 \int_\Omega u^2  - \beta_1 \int_\Omega v^2.
\end{equation}
Combining \eqref{Bag00} and \eqref{Bag01}, Young's inequality provides (recall $\alpha \chi\leq \beta_1 < \min \{b, 4c\}$):
\begin{equation*}
\begin{split}
\frac{d}{dt} \int_\Omega u  + \int_\Omega u & = -c \int_\Omega |\nabla u|^2  - \beta_1 \int_\Omega (\nabla u \cdot \nabla v) - \beta_1 \int_\Omega |\nabla v|^2 
- (b - \beta_1) \int_\Omega u^2  \\ & + (a + 1) \int_\Omega u  + h \int_{\partial \Omega} (\chi \alpha - \beta_1) uv 
- \beta_1 h \int_{\partial \Omega} v^2 - \beta_1 \int_\Omega v^2 
\\ &
\leq \int_\Omega \left[ -\left( \sqrt{c} \nabla u + \frac{\beta_1}{2\sqrt{c}} \nabla v \right) \cdot \left( \sqrt{c} \nabla u + \frac{\beta_1}{2\sqrt{c}} \nabla v \right) + \beta_1 \left( \frac{\beta_1}{4c} - 1 \right) |\nabla v|^2 \right] 
\\ & - \left(\frac{b - \beta_1}{2}\right) \int_\Omega u^2  + \const{KB} \leq \const{KB}
\quad \textrm{for all }  t \in (0,\TM).
\end{split}
\end{equation*}
Therefore, we have
\begin{equation*}\label{C0}
\frac{d}{dt} \int_\Omega u  + \int_\Omega u  \leq \const{KB} \quad \textrm{for all } t \in (0,\TM),
\end{equation*}
which implies the claim.
\end{proof}
\end{lemma}


\subsection{Uniform-in-time boundedness of the mass for $\tau=1$}
\begin{lemma}\label{BoundednessMassParab}
Let $\tau=1$, $\chi, a, h >0$, $\alpha \in [0,1]$, $C_{\partial \Omega}$ be the trace constant in Lemma \ref{LemmaTrace} and either 
\begin{equation}\label{AssumYuya}
b >
\begin{cases}
c + \frac{\alpha^2 \chi^2 h C_{\partial \Omega}}{16c} & \text{if } 0<c < \frac{\alpha \chi \sqrt{h C_{\partial \Omega}}}{4}, \\
\frac{\alpha \chi \sqrt{hC_{\partial \Omega}}}{2} & \text{if } c \geq \frac{\alpha \chi \sqrt{hC_{\partial \Omega}}}{4},
\end{cases}
\end{equation}
or
\begin{equation}\label{ConditionParaBagahei}
b > \chi.
\end{equation}
Then $u \in L^\infty((0,\TM);L^1(\Omega)).$ 
\begin{proof}
By reasoning as in the proof of Lemma \ref{BoundednessMass} and exploiting 
Lemma \ref{LemmaTrace} (with $\psi=u$, $\mathfrak{p}=2$), we arrive at
\[
\begin{split}
\frac{d}{dt} \int_{\Omega} u &\leq \frac{\alpha^2 \chi^2 h}{4 \epsilon_1} \int_{\partial \Omega} u^2 + \epsilon_1 h \int_{\partial \Omega} v^2 + a \int_{\Omega} u - b \int_{\Omega} u^2 - c \int_{\Omega} |\nabla u|^2 \\
& \leq \left(\frac{\alpha^2 \chi^2 h C_{\partial \Omega}}{4 \epsilon_1}-b\right) \int_{\Omega} u^2 
+ \epsilon_1 h \int_{\partial \Omega} v^2 + a \int_{\Omega} u +
\left(\frac{\alpha^2 \chi^2 h C_{\partial \Omega}}{4 \epsilon_1}- c\right) \int_{\Omega} 
|\nabla u|^2 \quad \textrm{for all }  t \in (0, \TM).
\end{split}
\]
Additionally, testing arguments employed on the the second equation in \eqref{problem} give thanks to Young's inequality 
\[
\frac{1}{2} \frac{d}{dt} \int_{\Omega} v^2 = -h \int_{\partial \Omega} v^2 - \int_{\Omega} |\nabla v|^2 - \int_{\Omega} v^2 + \int_{\Omega} u v
\leq -h \int_{\partial \Omega} v^2 - \int_{\Omega} v^2 + \frac{\epsilon_2}{\epsilon_1} \int_{\Omega} u^2 + \frac{\epsilon_1}{4 \epsilon_2} \int_{\Omega} v^2 \quad \textrm{on }   (0, \TM),
\]
so that both procedures provide
\begin{equation*}
\begin{split}
\frac{d}{dt} \int_{\Omega} \left( u + \frac{\epsilon_1}{2} v^2 \right) 
&\leq a \int_{\Omega} u + \left( \frac{\alpha^2 \chi^2 h C_{\partial \Omega}}{4 \epsilon_1} + \epsilon_2 - b \right) \int_{\Omega} u^2 
\\&
+ \left( \frac{\alpha^2 \chi^2 h C_{\partial \Omega}}{4 \epsilon_1} - c \right) \int_{\Omega} |\nabla u|^2 - \epsilon_1 \left( 1 - \frac{\epsilon_1}{4 \epsilon_2} \right) \int_{\Omega} v^2 \quad \textrm{on } (0,\TM).
\end{split}
\end{equation*}
At this stage, by taking into account \eqref{AssumYuya}, we can follow the proof of Lemma \ref{BoundednessMass} to choose $\epsilon_1$ and $\epsilon_2$ such that
\[
\frac{\chi^2 \alpha^2 h C_{\partial \Omega}}{4 \epsilon_1} + \epsilon_2 - b < 0, \quad \frac{\chi^2 \alpha^2 h C_{\partial \Omega}}{4 \epsilon_1} - c \leq 0 \quad \text{and} \quad 1 - \frac{\epsilon_1}{4 \epsilon_2} > 0.
\]
An application of Young's inequality allows to derive the following estimate 
\[
\frac{d}{dt} \int_{\Omega} \left( u + \frac{\epsilon_1}{2} v^2 \right) \leq -\const{Yy1} \int_{\Omega} \left( u + \frac{\epsilon_1}{2} v^2 \right) + \const{Yy2} \quad \text{on} \quad (0, \TM),
\]
which together with $\int_{\Omega} \left( u(x,0) + \frac{\epsilon_1}{2} v^2(x,0) \right)dx=\int_{\Omega} \left( u_0(x) + \frac{\epsilon_1}{2} v_0^2(x) \right)dx$ gives the conclusion.

We can achieve the claim also exploiting constrain \eqref{ConditionParaBagahei}.
For all $(x,t)\in \Omega \times (0,\TM)$, let us define the following function:
\[
\phi(v) = e^{mv} \quad \text{with} \quad \chi \leq m < b.
\]
We, subsequently, have
\[
\phi'(v) = m \phi(v) \quad \textrm{on} \; \Omega \times (0,\TM), 
\]
so that for $\mu \geq \frac{1}{c}$ 
it occurs that 
\[
\frac{d}{dt}  \int_\Omega (u + \mu) \phi(v)= \int_\Omega \phi(v) u_t  + \int_\Omega (u + \mu) \phi'(v) v_t  = \int_\Omega \phi(v) u_t + m \int_\Omega (u + \mu) \phi(v) v_t \quad \textrm{on }  (0,\TM).
\]
Using integration by parts, we obtain
\begin{equation}\label{bagRel1} 
\begin{split}
\frac{d}{dt}  \int_\Omega (u + \mu) \phi(v)+ \int_\Omega (u + \mu) \phi(v) & = \int_\Omega \left[ -c |\nabla u|^2 - 2m \, \nabla u \cdot \nabla v + m \left( u (\chi - m) - m \mu \right) |\nabla v|^2 \right] \phi(v)
\\ &
\quad + h \int_{\partial \Omega} v \left[ u (\alpha \chi - m) - m \mu \right] \phi(v)\\&
\quad +\int_{\Omega} \left[\left( - b + m \right) u^2 + \left(  (a + m\mu + 1) - mv \right) u \right]\phi(v)  + \mu \int_{\Omega} \left(1 - mv \right) \phi(v) \\ & := I_1 + I_2 + I_3 + I_4
\quad \textrm{for every } \; t<\TM.
\end{split}
\end{equation}
We now control $I_1$ and $I_2$. Using $\chi \leq m$ and $\mu \geq \frac{1}{c}$, we find on $(0,\TM)$
\begin{equation*}
\begin{split}
I_1 &= \int_{\Omega} \left[- c |\nabla u|^2 - 2m \, \nabla u \cdot \nabla v + m \left( u (\chi - m) - m\mu \right) |\nabla v|^2 \right] \phi(v) \\
&
= \int_{\Omega} \left[-\left(\sqrt{c} \nabla u + \frac{m}{\sqrt{c}} \nabla v \right) \cdot \left( \sqrt{c} \nabla u + \frac{m }{\sqrt{c}} \nabla v \right)
+ m \left( m \left (\frac{1}{c} - \mu \right) + u (\chi - m) \right) |\nabla v|^2 \right] \phi(v),
\end{split}
\end{equation*}
and
\[
I_2 = h \int_{\partial \Omega} v \left[ u (\alpha \chi - m) - m\mu \right] \phi(v),
\]
so that 
\begin{equation}\label{StimeI12}
I_1 \leq 0 \quad \textrm{and}\quad I_2\leq 0 \quad \textrm{on}  \; (0,\TM).
\end{equation}
Next, in order to control the term $I_3$, for each $t>0$, we divide $\Omega$ into two sets:  
\begin{equation*}
\Omega_1=\left\{x \in \Omega: v \geq \frac{a + m\mu + 1}{m} \right\} \quad \textrm{and} \quad \Omega_2=\left\{x \in \Omega: v < \frac{a + m\mu + 1}{m} \right\}.
\end{equation*}
We have
\begin{equation}\label{StimaI3}
\begin{split}
I_3 &= \int_{\Omega} \left[ (-b + m) u^2 + \left( (a + m\mu + 1) - mv \right) u \right] \phi(v) 
\\ &
= \int_{\Omega_1} \left[(-b + m) u^2 + \left( (a + m\mu + 1) - mv \right) u \right] \phi(v) 
+ \int_{\Omega_2} \left[ (-b + m) u^2 + \left( (a + m\mu + 1) - mv \right) u \right]\phi(v) 
\\ &
\leq \int_{\Omega_2} \left[(-b + m) u^2 + (a + m\mu + 1)u \right] \phi(v) 
\\ &
\leq \int_{\Omega_2 \cap \left\{ u \geq \frac{a + m\mu + 1}{b - m} \right\}} \left[ -\left( b - m \right) u + (a + m\mu + 1) \right] u \, \phi(v) 
 + e^{a+m\mu+1} \int_{\Omega_2 \cap \left\{ u < \frac{a + m\mu + 1}{b - m} \right\}} \left[ -\left( b - m \right) u + (a + m\mu + 1) \right] u 
\\ &
\leq \frac{e^{a+m\mu+1} (a + m\mu + 1)^2 |\Omega|}{b - m}\quad \textrm{for all } t \in (0,\TM).
\end{split}
\end{equation}
In order to control $I_4$, for each $t>0$ we similarly divide $\Omega$ into two sets:
\begin{equation*}
\Omega_3= \left\{x \in \Omega: v \geq \frac{1}{m} \right\}
\quad \textrm{and}\quad 
\Omega_4=\left\{x\in \Omega: v < \frac{1}{m} \right\}.
\end{equation*}
Subsequently, we obtain 
\begin{equation}\label{StimaI4}
I_4 = \mu \int_{\Omega} \left(1 - mv \right) \phi(v) 
= \mu \int_{\Omega_3} \left(1 - mv \right) \phi(v)  + \mu \int_{\Omega_4} \left(1 - mv \right) \phi(v) 
\leq \mu e |\Omega| \quad \textrm{for all } t \in (0,\TM).
\end{equation}
Considering estimates \eqref{StimeI12}, \eqref{StimaI3} and \eqref{StimaI4} for $I_i$ ($i=1,2,3,4$), relation \eqref{bagRel1} is turned into 
\[
\frac{d}{dt}  \int_\Omega (u + \mu) \phi(v)+ \int_\Omega (u + \mu) \phi(v) \leq \const{bag} \quad \textrm{for all } t\in (0,\TM),
\]
so that from $\phi(v)\geq 1$, we have the claim.
%
%
%
%
%
\end{proof}
\end{lemma}
\subsection{Uniform-in-time boundedness in $L^p(\Omega)$ for $\tau=0$}
Even though in the case of zero total flux ($\alpha=0$) the mass can be controlled even in the absence of classical logistic sources, their presence is not sufficient to provide the same estimate when a positive total flux ($\alpha>0$) is introduced. Indeed, as discussed in the $\S$\ref{SectionDampeningGrad}, positive fluxes tend to promote mass accumulation; to mitigate uncontrolled growth, it is essential to incorporate logistic sources including significant damping. The situation is even more evident when  boundedness of $L^p(\Omega)$-norms are required; in this sense, the parameter \( c \) associated with the gradient term in the subsequent lemmas (Lemmas  \ref{LemmaEll} and \ref{LemmaParab}) must be understood as strictly positive, even in the limit case \( \alpha = 0 \).

In the next lines, and for some $p>1$, we will make use of the functional $\frac{1}{p}\int_\Omega u^p$ for all $t \in (0,\TM)$.
\begin{lemma}\label{LemmaEll}
Let $n \in \mathbb{N}$ and the hypotheses of Lemma \ref{BoundednessMass} be complied. Then $u \in L^{\infty}((0,\TM); L^p(\Omega))$ for all $p>1$.
\begin{proof}
By multiplying the first equation in \eqref{problem} by $u^{p-1}$ and by integrating it over $\Omega$, we have
\begin{align}\label{ene:Lp:1}
\frac1p \frac{d}{dt} \int_\Omega u^p &= \int_{\partial \Omega} u^{p-1} (u_\nu - \chi u v_\nu) - (p-1) \int_\Omega u^{p-2}|\nabla u|^2
+ (p-1) \chi \int_\Omega u^{p-1} \nabla u \cdot \nabla v \notag \\
&\quad\, + a \int_\Omega u^p - b \int_\Omega u^{p+1} - c \int_\Omega u^{p-1}|\nabla u|^2 \notag \\
&= \alpha \chi h \int_{\partial \Omega} u^p v - (p-1) \int_\Omega u^{p-2}|\nabla u|^2 + \frac{p-1}p \chi \int_{\partial \Omega} u^p v_\nu - \frac{p-1}p \chi \int_\Omega u^p \Delta v
\notag \\
&\quad\,+ a \int_\Omega u^p - b \int_\Omega u^{p+1} - \frac{4c}{(p+1)^2} \int_\Omega |\nabla u^{\frac{p+1}2}|^2 \quad \textrm{for all } t \in (0,\TM).
\end{align}
As to the boundary integrals, recalling the Robin condition in \eqref{problem} we have
\begin{equation}\label{BoundaryNeg}
\int_{\partial \Omega} u^p v_\nu = - h \int_{\partial \Omega} u^p v \le 0 \quad \textrm{on $(0,\TM)$},
\end{equation}
so that the remaining integral term on the boundary of $\Omega$ to be treated is $\alpha \chi h  \int_{\partial \Omega} u^p v$.
A combination of Young's and trace inequalities (recall Lemma \ref{LemmaTrace} with $\psi=u^{\frac{p+1}{2}}$, $\psi=v^{\frac{p+1}{2}}$ and $\mathfrak{p}=2$) yields for all $t \in (0,\TM)$
\begin{equation}\label{BoundaryInt0}
\begin{split}
\alpha \chi h  \int_{\partial \Omega} u^p v &\leq \frac{2c}{(p+1)^2 C_{\partial \Omega}}
\int_{\partial \Omega} u^{p+1} + \const{Bb} \int_{\partial \Omega} v^{p+1} 
\\
  &  \leq \frac{2c}{(p+1)^2} \int_\Omega u^{p+1} + \frac{2c}{(p+1)^2}
  \int_\Omega |\nabla u^{\frac{p+1}{2}}|^2 + \const{BbF}\left( \int_{\Omega} v^{p+1}+\int_{\Omega}  |\nabla v^{\frac{p+1}{2}}|^2\right)\\
  &  
  \leq \frac{2c}{(p+1)^2} \int_\Omega u^{p+1} + \frac{2c}{(p+1)^2}
  \int_\Omega |\nabla u^{\frac{p+1}{2}}|^2 + \const{BbG}\left( \int_{\Omega} v^{p+1}+\int_{\Omega}  |\nabla v|^{p+1}\right)
  \\
  &  
  \leq \const{BbH} \int_\Omega u^{p+1} + \frac{2c}{(p+1)^2}
  \int_\Omega |\nabla u^{\frac{p+1}{2}}|^2,
\end{split}
\end{equation}
where we have applied in the last step relation \eqref{ellipticreg}, so being $\const{BbH}$ depending also on $C_E$.
Recalling that $\Delta v=v-u$ for all $(x,t) \in \Omega \times (0,\TM)$, by plugging bound \eqref{BoundaryInt0} into estimate \eqref{ene:Lp:1} and exploiting the Young inequality to absorb the term $a \int_\Omega u^p$ with $-b \int_\Omega u^{p+1}$, we have on $(0, \TM)$
\begin{equation}\label{estimate2}
\begin{split}
\frac1p \frac{d}{dt} \int_\Omega u^p & \leq 
 -\frac{4(p-1)}{p^2} \int_\Omega |\nabla u^{\frac{p}{2}}|^2 +\const{j} \int_\Omega u^{p+1} -\frac{2c}{(p+1)^2}\int_\Omega |\nabla u^{\frac{p+1}{2}}|^2 + \const{e}.
\end{split}
\end{equation}
Now we treat the term $\const{j} \int_\Omega u^{p+1}$ by exploiting a combination of the Gagliardo--Nirenberg and Young's inequalities. Setting $\theta_1:= \frac{p}{p+\frac{2}{n}} \in (0,1)$ and taking into account the boundedness of the mass given in Lemma \ref{BoundednessMass} yield
\begin{equation} \label{GN4}
\begin{split}
\const{j} \int_\Omega u^{p+1}&=\const{j}\|u^{\frac{p+1}{2}}\|_{L^2(\Omega)}^2
\leq \const{gn114} \|\nabla u^{\frac{p+1}{2}}\|_{L^2(\Omega)}^{2 \theta_1}  
\|u^{\frac{p+1}{2}}\|_{L^\frac{2}{p+1}(\Omega)}^{2(1-\theta_1)} + \const{gn114}   
\|u^{\frac{p+1}{2}}\|_{L^\frac{2}{p+1}(\Omega)}^2 \\
& \leq \const{gn124} \left(\int_\Omega |\nabla u^{\frac{p+1}{2}}|^2\right)^{\theta_1}+ \const{gn134}\leq \frac{2c}{(p+1)^2} \int_\Omega |\nabla u^{\frac{p+1}{2}}|^2 + \const{A4}
\quad \text{on }(0,\TM).
\end{split}
\end{equation}
By invoking bound \eqref{GN4}, the estimate \eqref{estimate2} becomes \begin{equation}\label{estimate3}
\frac1p \frac{d}{dt} \int_\Omega u^p \leq -\frac{4(p-1)}{p^2} \int_\Omega |\nabla u^{\frac{p}{2}}|^2  + \const{ff} \quad \textrm{for all } t \in 0,\TM).
\end{equation}
An application of the Gagliardo--Nirenberg and the Young inequalities yields for $\theta_2:= \frac{p-1}{p-1+\frac{2}{n}} \in (0,1)$ (and also recalling the boundedness of the mass) the following estimate 
\begin{equation}\label{GN2}
\begin{split}
\frac{1}{p}\int_\Omega u^p&= \frac{1}{p} \|u^{\frac{p}{2}}\|_{L^2(\Omega)}^2 \leq \const{gn21} \|\nabla u^{\frac{p}{2}}\|_{L^2(\Omega)}^{2 \theta_2}  
\|u^{\frac{p}{2}}\|_{L^\frac{2}{p}(\Omega)}^{2(1-\theta_2)} + \const{gn21} \|u^{\frac{p}{2}}\|_{L^\frac{2}{p}(\Omega)}^2 \\
& \leq \const{GN21} \left(\int_\Omega |\nabla u^{\frac{p}{2}}|^2\right)^{\theta_2}+ \const{Gn21} \leq \frac{4(p-1)}{p^2} \int_\Omega |\nabla u^{\frac{p}{2}}|^2 + \const{B} \quad \text{for all } t \in (0,\TM).
\end{split}
\end{equation}
By virtue of \eqref{GN2}, we derive from \eqref{estimate3} 
\begin{equation*}
\frac{d}{dt}\int_\Omega u^p \leq \const{l} -  \int_\Omega u^p  \text{ on }(0,\TM),
\end{equation*}
and henceforth $u \in L^{\infty}((0,\TM); L^p(\Omega))$ for all $p>1$.
\end{proof}
\end{lemma}
\subsection{Uniform-in-time boundedness in $L^p(\Omega)$ for $\tau=1$}
\begin{lemma}\label{LemmaParab}
Let $n \in \mathbb{N}$ and the hypotheses of Lemma \ref{BoundednessMassParab} be complied. Then $u \in L^{\infty}((0,\TM); L^p(\Omega))$ for all $p>1$.
\begin{proof}
By adhering to the line of reasoning presented for \eqref{ene:Lp:1} and also taking into account \eqref{BoundaryNeg}, since from the Young inequality 
it follows that 
\begin{equation*}\label{ene:Lp:2}
\alpha \chi h \int_{\partial \Omega} u^p v \le \frac{2c}{(p+1)^2 C_{\partial \Omega}} \int_{\partial \Omega} u^{p+1} + \mathcal{C} h \int_{\partial \Omega} v^{p+1} 
\quad\textrm{on $(0,\TM)$},
 \end{equation*}
being $\mathcal{C}:= (p+1)^{p-1} (\alpha \chi)^{p+1} \left(\frac{h p}{2c} C_{\partial \Omega}\right)^p$, we have adding $\frac{1}{p} \int_\Omega u^p$ to both sides 
\begin{align}\label{Parab1}
\frac1p \frac{d}{dt} \int_\Omega u^p + \frac{1}{p} \int_\Omega u^p& \leq \frac{2c}{(p+1)^2 C_{\partial \Omega}} \int_{\partial \Omega} u^{p+1} + \mathcal{C} h \int_{\partial \Omega} v^{p+1}
- (p-1) \int_\Omega u^{p-2}|\nabla u|^2 - \frac{p-1}p \chi \int_\Omega u^p \Delta v
\notag \\
&\quad\,+ \left(a+ \frac{1}{p}\right) \int_\Omega u^p - b \int_\Omega u^{p+1} - \frac{4c}{(p+1)^2} \int_\Omega |\nabla u^{\frac{p+1}2}|^2 \quad \textrm{for all } t \in (0,\TM).
\end{align}
Similar testing arguments for the second equation in \eqref{problem} and Young's inequality imply that
\begin{equation*}\label{Stima_v_P}
\begin{split}
\frac1{p+1} \frac{d}{dt} \int_\Omega v^{p+1} &= \int_{\partial \Omega} v^p v_\nu - p \int_\Omega v^{p-1} |\nabla v|^2- \int_\Omega v^{p+1}+ \int_\Omega u v^p \\ 
&\le - h \int_{\partial \Omega} v^{p+1}- \frac1{p+1}\int_\Omega v^{p+1}+ \frac1{p+1}\int_\Omega u^{p+1}
\quad\textrm{on $(0,\TM)$},
\end{split}
\end{equation*}
which together with \eqref{Parab1} and $\left(a+\frac{1}{p}\right) \int_\Omega u^p \leq b \int_\Omega u^{p+1} + \const{GGG1}$ yields
\begin{align}\label{ene:Lp:3}
 &\frac{d}{dt} \left( \frac1p \int_\Omega u^p + \frac{\mathcal{C}}{p+1} \int_\Omega v^{p+1} \right) + \frac{1}{p} \int_\Omega u^p \notag \\
 &\quad\, \le \frac{2c}{(p+1)^2 C_{\partial \Omega}} \int_{\partial \Omega} u^{p+1}
 - \frac{p-1}p \chi \int_\Omega u^p \Delta v - \frac{4c}{(p+1)^2} \int_\Omega |\nabla u^{\frac{p+1}2}|^2\notag \\
&\quad\,\quad\, - \frac{\mathcal{C}}{(p+1)} \int_\Omega v^{p+1}
+ \frac{\mathcal{C}}{(p+1)}  \int_\Omega u^{p+1} + \const{GGG1} \quad\textrm{for all $t \in (0,\TM)$}.
  \end{align}
From Lemma \ref{LemmaTrace} (with $\psi=u^{\frac{p+1}{2}}$ and $\mathfrak{p}=2$) and Young's inequality, we derive 
\begin{equation*}\label{Parab2}
\frac{2c}{(p+1)^2 C_{\partial \Omega}} \int_{\partial \Omega} u^{p+1} - \frac{p-1}p \chi \int_\Omega u^p \Delta v \leq
\const{GG1} \int_{\Omega} u^{p+1} + \frac{2c}{(p+1)^2} \int_{\Omega} |\nabla u^{\frac{p+1}2}|^2 + \const{GVc1}
\int_\Omega |\Delta v|^{p+1} \quad \textrm{on } (0,\TM),
\end{equation*}
which plugged in \eqref{ene:Lp:3} entails 
\begin{align*}
&\frac{d}{dt}\left( \frac1p \int_\Omega u^p + \frac{\mathcal{C}}{p+1} \int_\Omega v^{p+1}\right)+ \left( \frac1p \int_\Omega u^p 
+ \frac{\mathcal{C}}{p+1}\int_\Omega v^{p+1} \right)\notag \\
&\quad\, \le - \frac{2c}{(p+1)^2} \int_\Omega |\nabla u^{\frac{p+1}2}|^2
+ \const{SG1} \int_\Omega u^{p+1} + \const{GVc1} \int_\Omega |\Delta v|^{p+1}
+ \const{GGG1} \quad \textrm{for all } t \in (0,\TM),
\end{align*}
or equivalently 
\begin{align*}
\frac{d}{dt} \left( e^t 
\left(\frac1p \int_\Omega u^p + \frac{\mathcal{C}}{p+1} \int_\Omega v^{p+1} \right)\right)&= e^t \frac{d}{dt} \left( \frac1p \int_\Omega u^p + \frac{\mathcal{C}}{p+1} \int_\Omega v^{p+1} \right)+ e^t \left(\frac1p \int_\Omega u^p + \frac{\mathcal{C}}{p+1}\int_\Omega v^{p+1} \right)\notag \\
&\quad\, \le \frac{-2c e^t}{(p+1)^2} \int_\Omega |\nabla u^{\frac{p+1}2}|^2
+ e^t \const{SG1} \int_\Omega u^{p+1} + e^t \const{GVc1}\int_\Omega |\Delta v|^{p+1}
+ e^t \const{GGG1} \quad \textrm{on } (0,\TM).
\end{align*}
Now we integrate the previous bound between $0$ and $t$, obtaining
\begin{align*}
&e^t \left(\frac1p \int_\Omega u^p + \frac{\mathcal{C}}{(p+1)} \int_\Omega v^{p+1}
\right) \\
&\quad\, \le \const{SF} - \frac{2c}{(p+1)^2} \int^t_0 e^{s} \left(
\int_\Omega |\nabla u^{\frac{p+1}2}(\cdot,s)|^2\right)\,ds
+ \const{SG1} \int^t_0 e^{s} \left( \int_\Omega u^{p+1}(\cdot,s) \right)\,ds\\
&\quad\,\quad\, + \const{GVc1} \int^t_0 e^{s} \left(\int_\Omega |\Delta v(\cdot,s)|^{p+1}
\right)\,ds+ \const{GGG1} (e^{t}-1) \quad \textrm{for all } t\in(0,\TM).
\end{align*}
By means of \eqref{tau1} in Lemma \ref{regularity} (with $z=v$, $q=p+1$ and $\psi=u$) it follows that
\begin{align}\label{ene:Lp:5}
&e^{t} \left(\frac1p \int_\Omega u^p + \frac{\mathcal{C}}{(p+1)} \int_\Omega v^{p+1}
\right) \notag \\
&\quad\, \le \const{SF1} - \frac{2c}{(p+1)^2}
\int^t_0 e^{s} \left(\int_\Omega |\nabla u^{\frac{p+1}2}(\cdot,s)|^2\right)\,ds
+ \const{Chicca} \int^t_0 e^{s} \left(\int_\Omega u^{p+1}(\cdot,s)\right)\,ds
\notag \\
&\quad\,\quad\,+ \const{GGG1} (e^{t}-1) \quad \textrm{for all } t\in(0,\TM).
\end{align}
By relying on relation \eqref{GN4} we have up to the constants
\begin{equation*} \label{GNS4}
\const{Chicca} \int_\Omega u^{p+1} \leq \frac{2c}{(p+1)^2} \int_\Omega |\nabla u^{\frac{p+1}{2}}|^2 + \const{As4}
\quad \textrm{on }(0,\TM),
\end{equation*}
which inserted into \eqref{ene:Lp:5}, leads to
\[
e^{t} \left(\frac1p \int_\Omega u^p + \frac{\mathcal{C}}{(p+1)} \int_\Omega v^{p+1}\right) \le \const{SF1}+ \const{GGG01} (e^{t}-1) \quad \textrm{for all } t\in(0,\TM),
\]
so giving the conclusion.
\end{proof}
\end{lemma}

\subsection{Proof of Theorem \ref{theoremlocal}}
For $\tau=0$, from the assumptions \eqref{Estimate_b} or \eqref{Estimate_Bag0}, Lemma \ref{LemmaEll} implies
$u \in L^{\infty}((0,\TM); L^p(\Omega))$ for all $p>1$, so that using \eqref{tau0extension} we obtain $\nabla v \in L^{\infty}((0,\TM); L^\infty(\Omega))$. Subsequently, Lemma \ref{LemmaMoserType} and Lemma \ref{ExtensionLemma} give the claim.
For $\tau=1$, whenever \eqref{Estimate_b} or \eqref{Estimate_Bag1} are complied, we have the same conclusion by exploiting Lemma \ref{LemmaParab} and relation \eqref{tau1extension}. 
\qed
\appendix
\section{The case $a=b=c=\alpha=0$: Boundedness and finite-time blow-up}\label{Appendix}
\resetconstants{c}
To highlight how the positive flow makes the analysis of the considered model more complex, in this appendix we analyze problem \eqref{problem} in the absence of logistics ($a=b=c=0$) and in the case of zero-flux ($\alpha=0$). In particular, we will see how the dynamics of the corresponding solutions to these mechanisms is similar to the one where totally insulated domains are considered. We are particularly referring to bounded and blow-up solutions to models with classic homogeneous Neumann boundary conditions. 

In this sense, for $\chi, h>0$ we will refer to the following Keller--Segel model 
\begin{equation}\label{problemB}
\begin{cases}
u_t= \Delta u - \chi \nabla \cdot (u\nabla v) &{\rm in}\ \Omega \times (0, \TM),\\
\tau v_t=\Delta v-v+ u &{\rm in}\ \Omega \times (0, \TM),\\
u_\nu - \chi u v_\nu =0, \; v_\nu=-hv &{\rm on}\ \partial\Omega \times (0, \TM),\\
u(x, 0)=u_0(x), \tau v(x,0)= \tau v_0(x) &x \in \bar{\Omega},
\end{cases} 
\end{equation}
mentioning that for the limit case $h=0$ (h.N.b.c.) it is known that 
\begin{itemize}
\item if $n=1$, all solutions are global in time and bounded (\cite{Nagai_1995}  for $\tau=0$, \cite{OsYagUnidim} for $\tau=1$),
\item if $n\geq 2$, there are initial data providing blow-up solutions (\cite{Nagai_1995} for $\tau=0$, \cite{CS-2014,W} for $\tau=1$).
\end{itemize}
In consideration of these factors, the question is raised as to whether the results of boundedness and blow-up remain applicable in the context of Robin boundary conditions; this is the objective of this section.

For our purpose, herein 
\begin{equation}\label{AssumpAnnex}
\begin{cases}
\Omega=B_R \;\textrm{is a ball in}\; \mathbb{R}^n, n\geq 1, R>0, \delta \in (0,1), \tau\in\{0,1\},\\  
 u_0, \tau v_0: \bar{\Omega}  \rightarrow \R^+ , \textrm{radially symmetric, with } u_0,  \tau v_0 \in C^{2+\delta}(\bar\Omega) \\
 \textrm{complying with} \quad u_{0\nu}=-\chi h u_0v_0 \quad \textrm{ and } \quad \tau v_{0\nu}=-h \tau v_0 \quad  \textrm{ on }\partial \Omega.
 \end{cases}
\end{equation}
By naturally exploiting the local existence result given in Lemma \ref{theoremExistence}, we obtain a couple of radially symmetric functions $(u,v)$ defined in $\Omega\times (0,\TM)$ solving classically problem \eqref{problemB}, where $\TM$ is such that 
\begin{equation}\label{ExtensionAnnexo}
\textrm{either}\,\; \TM=\infty \quad \textrm{or} \quad \limsup_{t \nearrow \TM}\lVert u(\cdot,t)\rVert_{L^\infty(\Omega)} =\infty. 
\end{equation}
Additionally, $u$ manifestly fulfills
\begin{equation}\label{MassZero}
\int_\Omega u(x,t)\,dx=\int_\Omega u_0(x)\,dx \quad \textrm{for all } t \in (0,\TM),
\end{equation}
so that by integrating the second equation of \eqref{problemB} and exploiting \eqref{MassZero}, this leads to
\begin{equation*}
\begin{split}
\frac{d}{dt} \int_\Omega \tau v(x,t)\,dx&=\int_{\partial \Omega} v_\nu-\int_{\Omega} v + \int_\Omega u = -\int_{\partial \Omega} h v -\int_{\Omega} v + \int_\Omega u_0(x)\,dx
 \leq  -\int_{\Omega} v + \int_\Omega u_0(x)\,dx \quad \textrm{on } (0,\TM),
\end{split}
\end{equation*}
which implies
\begin{equation}\label{VZero}
\int_\Omega v(x,t)\, dx \leq \max\left\{\int_\Omega u_0(x)\,dx, \int_\Omega \tau v_0(x)\,dx\right\} \quad \textrm{for all } t \in (0,\TM).
\end{equation}
\subsection{Boundedness in the one-dimensional case}
\begin{theorem}\label{theoremlocalN1}
Let the hypotheses in \eqref{AssumpAnnex} be fulfilled, $n=1$ and $\chi, h >0$. Then 
problem \eqref{problemB} admits a unique solution
\begin{equation*}
(u,v)\in C^{2+\delta,1+\frac{\delta}{2}}( \Bar{\Omega} \times [0, \infty))\times C^{2+\delta,\tau+\frac{\delta}{2}}( \Bar{\Omega} \times [0, \infty))
\end{equation*}
such that $0\leq u,v \in L^\infty(\Omega \times (0,\infty)).$ 
\begin{proof}
From \eqref{ExtensionAnnexo} and standard regularity results it is sufficient to show that $u \in L^{\infty}((0,\TM); L^p(\Omega))$ for all $p>1$. We distinguish the cases $\tau=0$ and $\tau=1$.

\textbf{Case $\tau=0$.} By taking $\alpha=a=b=c=0$ in relation \eqref{ene:Lp:1} and recalling that $\Delta v= v -u$ from the second equation of problem \eqref{problemB}, we have
\begin{equation*}
\begin{split}
\frac1p \frac{d}{dt} \int_\Omega u^p &= - (p-1) \int_\Omega u^{p-2}|\nabla u|^2 - \frac{p-1}p \chi \int_\Omega u^p \Delta v \leq - \frac{4(p-1)}{p^2} \int_\Omega 
|\nabla u^{\frac{p}{2}}|^2 + \frac{p-1}p \chi \int_\Omega u^{p+1} \quad \textrm{for all } t \in (0,\TM).
\end{split}
\end{equation*}
Set $\theta_3:= \frac{\frac{p}{2}-\frac{p}{2(p+1)}}{\frac{p}{2}+\frac{1}{2}} \in (0,1)$.
By virtue of the Gagliardo--Nirenberg and Young's inequalities and the boundedness of the mass \eqref{MassZero}, one obtains 
\begin{equation}\label{GNlast}
\begin{split}
\frac{p-1}p \chi \int_\Omega u^{p+1}&=\frac{p-1}p \chi \|u^{\frac{p}{2}}\|_{L^{\frac{2(p+1)}{p}}(\Omega)}^{\frac{2(p+1)}{p}}
\leq\const{ggSS1} \|\nabla u^{\frac{p}{2}}\|_{L^2(\Omega)}^{\frac{2(p+1)}{p} \theta_3}  
\|u^{\frac{p}{2}}\|_{L^{\frac{2}{p}}(\Omega)}^{\frac{2(p+1)}{p}(1-\theta_3)} + \const{ggSS1} \|u^{\frac{p}{2}}\|_{L^\frac{2}{p}(\Omega)}^{\frac{2(p+1)}{p}}\\ 
&\leq \const{ggSS2} \left(\int_\Omega |\nabla u^{\frac{p}{2}}|^2\right)^{\frac{(p+1)}{p}\theta_3}+ \const{ggS1} \leq \frac{(p-1)}{p^2} \int_\Omega |\nabla u^{\frac{p}{2}}|^2
+ \const{ggS12} \quad \textrm{on } (0,\TM), 
\end{split}
\end{equation}
which, inserted in the above inequality, gives
\begin{equation*}
\begin{split}
\frac1p \frac{d}{dt} \int_\Omega u^p \leq - \frac{3(p-1)}{p^2} \int_\Omega 
|\nabla u^{\frac{p}{2}}|^2 + \const{ggS12} \quad \textrm{for all } t \in (0,\TM).
\end{split}
\end{equation*}
At this point we obtain the claim by following the proof of Lemma \ref{LemmaEll}.

\textbf{Case $\tau=1$.} Similarly to the previous case, by following the proof of Lemma \ref{LemmaParab} with $\alpha=a=b=c=0$, we arrive at 
\begin{equation*}
\begin{split}
e^{t} \left(\frac1p \int_\Omega u^p + \frac{1}{(p+1)} \int_\Omega v^{p+1}
\right) &\le \const{SFgg1} - \frac{4(p-1)}{p^2}
\int^t_0 e^{s} \left(\int_\Omega |\nabla u^{\frac{p}{2}}(\cdot,s)|^2\right)\,ds
+ \const{Chicca11} \int^t_0 e^{s} \left(\int_\Omega u^{p+1}(\cdot,s)\right)\,ds\\
&+ \const{GGGss1} (e^{t}-1) \quad \textrm{on } (0,\TM).
\end{split}
\end{equation*}
By plugging bound \eqref{GNlast} (up to the constants) in the previous inequality, we derive 
\[
e^{t} \left(\frac1p \int_\Omega u^p + \frac{1}{(p+1)} \int_\Omega v^{p+1}\right) \le \const{SFgg1} + \const{GGG1c} (e^{t}-1) \quad \textrm{for all } t\in(0,\TM),
\]
so giving the statement.
\end{proof}
\end{theorem}
\subsection{Finite-time blow-up}
In this section we will assume $n\geq 2$.
\subsubsection{Blow-up for $\tau=0$}
The following derivations are obtained by adapting arguments in \cite{Nagai_1995}, where homogeneous Neumann boundary conditions for both $u$ and $v$ are set.

Since $u$ and $v$ are radially symmetric, we can write $u(r,t)$ and $v(r,t)$ with $r=|x|$ instead of $u(x,t)$ and $v(x,t)$, respectively. Based on the same \cite{Nagai_1995}, we define the functions $U$ and $V$ as 
\begin{align}\label{def:UV}
U(r,t) := \int^r_0 \rho^{n-1} u(\rho,t) \,d\rho
\quad\textrm{and} \quad V(r,t) := \int^r_0 \rho^{n-1} v(\rho,t) \,d\rho
\end{align}
for $r\in[0,R]$ and $t\in[0,\TM)$. Moreover, being $\omega_n$ the area of the unit sphere in $\mathbb{R}^n$, we will make mention to the following quantities:  
\[
M_n(t) := \frac1{\omega_n} \int_\Omega u(x,t)|x|^n \,dx
\quad\textrm{for $t\in[0,\TM)$}
\]
and 
\[
\theta := \frac{1}{\omega_n} \int_\Omega u_0(x) \,dx. 
\]
\begin{lemma}\label{ddtM_N}
The following inequality holds:
\[
\frac{d}{dt} M_n(t) \le 2n(n-1) \theta^{\frac2n} M_n^{\frac{n-2}n}(t)
- \frac{n}2 \chi \theta^2 + \chi n \int^R_0 r^{n-1} u(r,t)V(r,t) \,dr \quad  \textrm{for all } t\in[0,\TM).
 \]
\begin{proof}
The boundary condition $u_\nu - \chi uv_\nu = 0$ 
ensures that \cite[(3.3)]{Nagai_1995} holds. Moreover, from \eqref{MassZero} we can derive $U(R,t) = \theta$ for all $t\in[0,\TM)$.
Thus this lemma can be proved by an argument similar to that of the proof of \cite[Lemma 3.2]{Nagai_1995}.
\end{proof}
\end{lemma}
\begin{lemma}\label{intuV}
The following inequality holds:
\[
n \chi \int^R_0 r^{n-1} u(r,t)V(r,t) \,dr \le \chi n R^{-n} \theta M_n(t) + \chi J_\theta(M_n(t)) \quad  \textrm{for all } t\in(0,\TM),
\]
where for $s\ge0$,
\begin{align*}
J_\theta(s):=
\begin{cases}
\dfrac1e \theta^\frac32 s^\frac12 &\textrm{if $n=2$},\\[3mm]
\dfrac{n}{2(n-2)} \theta^\frac{2n-2}{n} s^\frac2n &\textrm{if $n\ge3$}.
\end{cases}
\end{align*}
\begin{proof}
As in the proof of \cite[Lemma 3.3]{Nagai_1995}, we put 
\[
\Phi(r,t) := V(r,t) - \theta \left( \frac{r}{R} \right)^n
\quad\textrm{for $r\in[0,R]$ and $t\in(0,\TM)$}.
\]
Moreover we define $w$ as 
\begin{align*}
w(r,t) := 
\begin{cases}
 - \dfrac{\theta}{2} r^2 \log{\dfrac{r}{R}} - h V_r(R,t) \left( \dfrac{r}{R} \right)^2
 &\textrm{if $n=2$},\\[3mm]
\dfrac{\theta R^{-n}}{2(n-2)} (R^n r^2 - R^2 r^n) - h V_r(R,t) \left( \dfrac{r}{R} \right)^2
&\textrm{if $n\ge3$}
\end{cases}
\end{align*}
for $r\in[0,R]$ and $t\in(0,\TM)$.  Let us fix $t\in(0,\TM)$. It follows from \eqref{def:UV} that $V_r = r^{n-1} v$ and 
$V_{rr} - \frac{n-1}{r}V_r = r^{n-1}v_r$ for all $r\in(0,R)$.
Also, the second equation in \eqref{problemB} is rewritten as 
\begin{align}\label{2nd:rad}
0 = (r^{n-1}v_r)_r - r^{n-1}v + r^{n-1}u
\quad\textrm{for all $r\in(0,R)$}.
\end{align} 
We integrate it to see that
$V_{rr} - \frac{n-1}{r} V_r - V = - U$ for all $r\in(0,R)$, which yields
\[
\Phi_{rr} - \dfrac{n-1}{r} \Phi_r - \Phi = \theta \left( \dfrac{r}{R} \right)^n - U
\quad \textrm{for all $r\in(0,R)$}.
\]
Moreover, integrating \eqref{2nd:rad} again, we infer from the equality $U(R,t) = \theta$ and the boundary condition $v_\nu = -hv$ in \eqref{problemB} 
that $V(R,t) - \theta = R^{n-1}v_r(R,t) = - h V_r(R,t)$.
Hence, the function $\Phi$ fulfills that 
\[
\Phi(0,t) = 0 \quad \textrm{and}\quad \Phi(R,t) = - h V_r(R,t).
\] 
On the other hand, from straightforward computations we see that the function $w$ satisfies
\begin{align*}
\begin{cases}
w_{rr} - \dfrac{n-1}{r} w_r - w = - \theta - w, \quad r\in(0,R),\\[3mm]
w(0,t) = 0, \quad w(R,t) = - h V_r(R,t). 
\end{cases}
\end{align*}
Therefore, from the comparison argument we obtain $\Phi(r,t) \le w(r,t)$ for all 
$r\in(0,R)$ and $t\in(0,\TM)$, that is, 
$V(r,t) \le \theta \left( \frac{r}{R} \right)^n + w(r,t)$ for all $r\in(0,R)$ and $t\in(0,\TM)$. 
Thus
\[
\chi n \int^R_0 r^{n-1} u(r,t)V(r,t) \,dr \le \chi n R^{-n} \theta M_n(t)
+ \chi n \int^R_0 r^{n-1} u(r,t)w(r,t) \,dr \quad \textrm{for all } t\in(0,\TM),
\]
which together with the estimate
\[
w(r,t) \le
\begin{cases}
\dfrac{\theta R}{2e}r &\textrm{if $n=2$},\\[3mm]
\dfrac{\theta}{2(n-2)}r^2 &\textrm{if $n\ge3$}
\end{cases}
\]
for all $r\in(0,R)$ and $t\in(0,\TM)$ and H\"{o}lder's inequality implies the claim.
\end{proof}
\end{lemma}
By the same procedure in \cite{Nagai_1995} we can obtain the result on finite-time blow-up. 
\begin{theorem}[Finite-time blow-up $\tau=0$]\label{thm2}
For $n \ge 2$, let the hypotheses in \eqref{AssumpAnnex} be complied. Moreover, assume that $\int_\Omega u_0(x)\, dx > \frac{8\pi}{\chi}$ when $n = 2$. Then there exists $c(\theta)>0$ such that 
if $0<\frac{1}{\omega_n} \int_\Omega u_0(x)|x|^n\, dx < c(\theta)$, 
then the solution $(u,v)$ of \eqref{problemB} blows up in finite time.
\begin{proof}
In order to prove finite-time blow-up by a contradiction argument, we assume $\TM = \infty$. Lemmas \ref{ddtM_N} and \ref{intuV} yield 
\begin{align}\label{ddtM_N:2}
\frac{d}{dt} M_n(t) \le E_\theta(M_n(t)) \quad\textrm{for all $t\in(0,\infty)$},
\end{align} 
where 
\[
E_\theta(s):= 2n(n-1) \theta^{\frac2n} s^{\frac{n-2}n} - \frac{n}2 \chi \theta^2
+ \chi n R^{-n} \theta s + \chi J_\theta(s) \quad\textrm{for $s\in(0,\infty)$}.
\]
Here, 
\[
E_\theta(0) =
\begin{cases}
(4 - \chi \theta) \theta &\textrm{if $n=2$}, \\
- \frac{n}2 \chi \theta^2 &\textrm{if $n\ge3$},
\end{cases}
\]
which implies that $E_\theta(0) < 0$ in the case that $n=2$ and $\theta > \frac{4}{\chi}$ and in the case that $n\ge3$ and $\theta>0$. 
Therefore we can find $c_1 = c_1(\theta)>0$ such that if $0< M_n(0) < c_1$, then $E_\theta(M_n(0))<0$. This together with \eqref{ddtM_N:2} implies that if $0< M_n(0) < c_1$, then there exists $T_0 \in (0,\infty)$ such that $M_n(t) \to 0$ as $t \to T_0$. 
Thus we see that $\TM \le T_0$, which contradicts with $\TM = \infty$. Hence from \eqref{ExtensionAnnexo} we arrive at the conclusion.
\end{proof}
\end{theorem}
\subsubsection{Blow-up for $\tau=1$}
For the analysis of the  fully parabolic scenario, we will adapt arguments in \cite{CS-2014} for $n=2$ and in \cite{W} for $n\geq3$, where homogeneous Neumann boundary conditions for both $u$ and $v$ are considered.  Since most necessary adjustments are rather small, we will confine ourselves to a brief outline.

Let us introduce for arbitrary smooth positive functions $u$ and $v$ and $\chi, h>0$, the Lyapunov  functional 
\begin{equation}\label{F}
\mathcal{F}_h(u,v):=\frac{\chi}{2}\int_\Omega |\nabla v|^2+\frac{\chi}{2}\int_\Omega v^2 - \chi \int_\Omega uv + \int_\Omega u \, \ln u +\frac{\chi h}{2}\int_{\partial \Omega} v^2
\end{equation}
and the dissipation rate  
\begin{equation}\label{D}
\mathcal{D}(u,v):= \chi \int_\Omega v_t^2 + \int_\Omega u \left|\frac{\nabla u}{u}-\chi \nabla v\right|^2.
\end{equation}
It can be noted that, conversely to the classical definition of the Lyapunov functional $\mathcal{F}$
(see for instance \cite[(1.5)]{W}, for the special case $\chi=1$), in our scenario an extra term has to be considered in the expression of the functional 
$\mathcal{F}_h$ itself; it is also seen that when $h \to 0$ and $\chi \to 1$ we have that $\mathcal{F}_h \to \mathcal{F}$.


Since $(u,v)$ is radially symmetric we deduce a pointwise upper bound for $v$.
\begin{lemma}\label{UpperBoundV}
Let $\kappa > n-2$ if $n \geq 3$, and $\kappa=2$ if $n=2$. Then one can find $C(\kappa)>0$ such that for all radially symmetric and positive functions $u_0, v_0 \in C^{2+\delta}(\bar{\Omega})$, the corresponding solution of \eqref{problemB} satisfies
\begin{equation*}
v(r,t)\leq C(\kappa) \left(\|u_0\|_{L^1(\Omega)}+\|v_0\|_{L^1(\Omega)}+
\|\nabla v_0\|_{L^2(\Omega)}\right) r^{-\kappa} \quad \textrm{for all } (r,t) \in (0,R) \times (0,\TM).
\end{equation*}
\begin{proof}
For $n\geq 3$, the proof is obtained by adapting reasoning developed in \cite[Lemma 3.1, Lemma 3.2 and Corollary 3.3]{W}; the same adaptations permit to show the result when $n=2$. To be more precise, the manipulations have to be performed in \cite[Lemma 3.1]{W}, which is derived by invoking properties of Neumann heat semigroup; in our context of Robin conditions, relations \eqref{SemigroupRobin2} and 
\eqref{SemigroupRobin3} play the same roles as in \cite[(3.2), (3.3)]{W}.
\end{proof}
\end{lemma}
One of the main tools to prove blow-up of the solutions of \eqref{problemB} consists in the following energy inequality, which is well known for the Neumann boundary conditions (see \cite[Lemma 2.1]{WinVolumeFill}).
\begin{lemma}\label{DerivativeEstimate}
Let $\TM$ be given in Lemma \ref{theoremExistence}. If $(u,v)$ is a classical solution of problem \eqref{problemB}, then for $\mathcal{F}_h(u,v)$ and $\mathcal{D}(u,v)$ defined in \eqref{F} and \eqref{D} respectively, this relation holds:
\begin{equation}\label{DerEst1}
\frac{d}{dt} \mathcal{F}_h(u(\cdot,t), v(\cdot,t)) = - \mathcal{D}(u(\cdot,t), v(\cdot,t))\quad \textrm{for all }  t \in (0, \TM).
\end{equation}
\begin{proof}
Let us compute the derivative with respect the time to the energy functional $\mathcal{F}_h(u,v)$. Since the pair $(u,v)$ solves \eqref{problemB}, we obtain for all $t \in (0, \TM)$
\begin{equation*}
\begin{split}
\frac{d}{dt} \mathcal{F}_h(u(\cdot,t), v(\cdot,t))&= -\chi \int_\Omega v_t \Delta v +
\chi \int_{\partial \Omega} v_t v_{\nu} +\chi \int_\Omega v v_t - \chi \int_\Omega v u_t
-\chi \int_\Omega u v_t + \int_\Omega u_t \, \ln u + \chi h \int_{\partial \Omega} v v_t\\
&= \chi \int_\Omega v_t (-\Delta v + v-u) + \int_\Omega u_t (\ln u- \chi v)\\
&= - \chi \int_\Omega v^2_t - \int_\Omega u |\nabla \ln u - \chi \nabla v|^2 + \int_{\partial \Omega} (\ln u-\chi v) (u_{\nu}-\chi u v_{\nu})= - \mathcal{D}(u,v).
\end{split}
\end{equation*}
\end{proof}
\end{lemma}
Throughout this section, let $m, M, B, \chi, h>0$ and $\kappa > n-2$ if $n\geq 3$ or $\kappa=2$ if $n=2$; let us also define the set
\begin{equation*}
\begin{split}
\mathcal{S}(m,M,B, \kappa, h):=&\left\{(u, v) \in C^1(\bar{\Omega}) \times C^2(\bar{\Omega}): u \textrm{ and } v \textrm{ are positive and radially symmetric}
\right.\\ 
& \left. \quad \textrm{with } v_\nu=-h v \textrm{ on } \partial \Omega 
 \textrm{ and such that } \int_\Omega u=m, \int_\Omega v \leq M \textrm{ and }\right.\\
& \left. \quad v(x) \leq B |x|^{-\kappa} \textrm{ for all } x \in \Omega\right\}.
\end{split}
\end{equation*}
Moreover, for $(u,v) \in \mathcal{S}(m,M,B, \kappa, h)$ we introduce the functions $f$ and $g$ 
\begin{equation}\label{fg}
f:= -\Delta v+v-u \quad \textrm{and} \quad g:= \left(\frac{\nabla u}{\sqrt{u}}-\chi \sqrt{u} \nabla v\right) \cdot \frac{x}{|x|} \quad (\textrm{for } x \neq 0),
\end{equation}
that can be also written, since $(u,v)$ is radial, as $f=-r^{1-n}(r^{n-1} v_r)_r+v-u$
and $g=\frac{u_r}{\sqrt{u}}-\chi \sqrt{u} v_r$. 
Now, we dedicate ourselves to estimate the functional $\mathcal{F}_h$ in terms of the dissipation rate $\mathcal{D}$ on the set $\mathcal{S}$. As the first step, we estimate
the integral $\int_\Omega uv$ in terms of $\mathcal{D}$: this is the goal of the next three lemmas.
\begin{lemma}\label{FirstEstimateInt_uv}
There exists $C(M,h)>0$ such that for all $(u,v) \in \mathcal{S}(m,M,B, \kappa, h)$ we have
\begin{equation}\label{FirstEstimateInt_uvF}
\int_\Omega uv \leq 2 \int_\Omega |\nabla v|^2 + C(M,h) \left(\|f\|_{L^2(\Omega)}^{\frac{2n+4}{n+4}}+1\right).
\end{equation}
\begin{proof}
By multiplying the expression of $f$ in \eqref{fg} for $v$ and by integrating over $\Omega$, we have
\begin{equation}\label{FirstEst1}
\int_\Omega uv = \int_\Omega |\nabla v|^2 + h \int_{\partial \Omega} v^2
+ \int_\Omega v^2 - \int_\Omega fv.   
\end{equation}
Since $\Omega$ is a ball, we can invoke \cite[Lemma A.1]{PAYNE-DisBoundaryIntergal2010} to control the integral term on the boundary, yielding after an application of the Young inequality
\begin{equation*}
h \int_{\partial \Omega} v^2 \leq \frac{h n}{R} \int_{\Omega} v^2 + 2 h \int_\Omega v |\nabla v| \leq \left(\frac{h n}{R} + 3 h^2 \right) \int_{\Omega} v^2 + \frac{1}{3} \int_\Omega |\nabla v|^2,
\end{equation*}
which, inserted in equality \eqref{FirstEst1}, gives
\begin{equation}\label{FirstEst2}
\int_\Omega uv \leq \frac{4}{3} \int_\Omega |\nabla v|^2 + \left(1+ \frac{hn}{R} + 3 h^2\right) \int_{\Omega} v^2 - \int_\Omega fv.   
\end{equation}
At this point we manage the terms involving $\int_{\Omega} v^2$ and $- \int_\Omega fv$ as in the proof of \cite[Lemma 4.2]{W}; this leads for some 
$\const{b2}=\const{b2}(M)$, $\const{b3}=\const{b3}(M)$ and 
$\const{b4}=\const{b4}(M)$ to
\begin{equation*}
\left(1+ \frac{h n}{R} + 3 h^2\right) \int_{\Omega} v^2 \leq \frac{1}{3} \int_\Omega |\nabla v|^2 + \const{b2} \left(1+ \frac{h n}{R} + 3 h^2\right)^{\frac{n+2}{2}}
\end{equation*}
and 
\begin{equation*}
-\int_{\Omega} f v \leq \frac{1}{3} \int_\Omega |\nabla v|^2 + \const{b3}\|f\|_{L^2(\Omega)}^{\frac{2n+4}{n+4}}+\const{b4}.
\end{equation*}
By plugging these two estimates into bound \eqref{FirstEst2}, we find
\begin{equation}\label{FirstEst3}
\int_\Omega uv \leq 2 \int_\Omega |\nabla v|^2 +  \max\left\{\const{b3},\const{b4} + \const{b2} \left(1+ \frac{h n}{R} + 3 h^2\right)^{\frac{n+2}{2}}\right\}(\|f\|_{L^2(\Omega)}^{\frac{2n+4}{n+4}}+1), 
\end{equation}
so concluding. 
\end{proof}
\end{lemma}
\begin{remark}[Comparison with the Neumann boundary conditions]
We point out that the presence of the positive constant $h$ in \eqref{FirstEst3} defining $C(M,h)$ in Lemma \ref{FirstEstimateInt_uv} is manifestly connected to the boundary condition under investigation. In particular, if $h \to 0$ we have that $\const{b2} \left(1+ \frac{h n}{R} + 3 h^2\right)^{\frac{n+2}{2}} \to \const{b2}$, the constant $C(M,h) \to C(M)$, so recovering $C(M)$ of \cite[Lemma 4.2]{W}. 
\end{remark}
In the next lemma, we will control the term  $\int_\Omega |\nabla v|^2$.
\begin{lemma}\label{EstimateGradOmega}
There exists $C(m, M, B, \kappa,\chi)>0$ such that each $(u,v) \in \mathcal{S}(m,M,B, \kappa, h)$ satisfies
\begin{equation}\label{EstimateGradOmegaF}
\int_{\Omega} |\nabla v|^2 \leq \frac{1}{4} \int_\Omega uv + 
C(m, M, B, \kappa,\chi) \left(\|f\|_{L^2(\Omega)}^{2\theta} +
\|g\|_{L^2(\Omega)}+ 1\right),
\end{equation}
where 
\begin{equation}\label{theta}
\theta:= \frac{1}{1+\frac{n}{(2n+4)\kappa}} \in \left(\frac{1}{2},1\right).
\end{equation}
\begin{proof}
For $n\geq 3$, and accordingly to \cite[Lemma 4.3 and Lemma 4.4]{W}, by splitting the domain $\Omega$ as $\Omega=B_{r_0} \cup (\Omega\setminus B_{r_0})$, where $B_{r_0}$ is a small inner ball, the estimate on $\int_{\Omega\setminus B_{r_0}}|\nabla v|^2$ is obtained by replacing \cite[(4.14)]{W} with 
\begin{equation*}
\alpha \int_\Omega v^{\alpha-1} |\nabla v|^2 = \int_\Omega f v^{\alpha} + 
\int_\Omega u v^{\alpha} - \int_\Omega v^{\alpha+1} - h \int_{\partial \Omega}  v^{\alpha+1} 
\leq  \int_\Omega f v^{\alpha} + \int_\Omega u v^{\alpha}. 
\end{equation*}
(For $n=2$ we refer to \cite[Lemma 2.3]{CS-2014}.)

On the other hand, $\int_{B_{r_0}} |\nabla v|^2$ is similarly controlled as in \cite[Lemma 4.4]{W}, for $n\geq 3$, and in \cite[Lemma 2.4]{CS-2014} for $n=2$.

From $\int_{\Omega} |\nabla v|^2= \int_{\Omega\setminus B_{r_0}} |\nabla v|^2 + \int_{B_{r_0}} |\nabla v|^2$ and following \cite[Lemma 4.5 (ArXiv version)]{W}, we definitively can conclude. 
%
%
%
%
%
\end{proof}
\end{lemma}
Now we are in the position to obtain the desired estimate of the term $\int_\Omega uv$.
\begin{lemma}\label{FinalEstimateInt_uv}
There exists $C(m,M,B, \kappa, h,\chi)>0$ such that for all $(u,v) \in \mathcal{S}(m,M,B, \kappa, h)$ we have
\begin{equation}\label{FinalEstimateInt_uvF}
\int_\Omega uv \leq C(m,M,B, \kappa, h,\chi) \left(\|f\|_{L^2(\Omega)}^{2\theta} +
\|g\|_{L^2(\Omega)}+1\right)
\end{equation}
with $\theta \in \left(\frac{1}{2},1\right)$ given in \eqref{theta}.
\begin{proof}
By exploiting bound \eqref{EstimateGradOmegaF} in estimate \eqref{FirstEstimateInt_uvF}, we obtain the claim, after an application of the Young inequality, since $\frac{2n+4}{n+4} \leq 2\theta$, exactly as in \cite[Lemma 4.1]{W}.
\end{proof}
\end{lemma}
Finally, the previous lemma allows us to estimate the Lyapunov functional $\mathcal{F}_h$ in terms of the dissipation rate $\mathcal{D}$ on the set $\mathcal{S}$, which is the other main tool to prove blow-up.
\begin{lemma}\label{TheoremBU}
There exists $C(m,M,B, \kappa, h,\chi)>0$ such that for all $(u,v) \in \mathcal{S}(m,M,B, \kappa, h)$ we have
\begin{equation*}\label{EstFD}
\mathcal{F}_h(u,v) \geq -C(m,M,B, \kappa, h,\chi) \, (\mathcal{D}^{\theta}(u,v)+1)
\end{equation*}
with $\theta \in \left(\frac{1}{2},1\right)$ given in \eqref{theta}.
\begin{proof}
By exploiting that $\xi \ln \xi \geq -\frac{1}{e}$ for all $\xi>0$ and bound \eqref{FinalEstimateInt_uvF}, we derive for some $\const{A1}=\const{A1}(m,M,B, \kappa, h,\chi)$
\begin{equation*}
\begin{split}
\mathcal{F}_h(u,v)&=\frac{\chi}{2}\int_\Omega |\nabla v|^2+\frac{\chi}{2}\int_\Omega v^2 - \chi \int_\Omega uv + \int_\Omega u \, \ln u +\frac{\chi h}{2}\int_{\partial \Omega} v^2\\
&\geq - \chi \int_\Omega uv  - \frac{|\Omega|}{e}\\
&\geq -\const{A1} \left(\|f\|_{L^2(\Omega)}^{2\theta} +
\|g\|_{L^2(\Omega)}+1\right).
\end{split}
\end{equation*}
Since $\theta>\frac{1}{2}$, an application of the Young inequality yields for some 
$\const{A2}=\const{A2}(m,M,B, \kappa, h,\chi)$ 
\begin{equation*}
\mathcal{F}_h(u,v) \geq -\const{A2} \left((\chi \|f\|_{L^2(\Omega)}^2 +
\|g\|_{L^2(\Omega)}^2)^{\theta}+1\right) 
= -\const{A2} \left(\mathcal{D}^{\theta}(u,v)+1\right),
\end{equation*}
where, by virtue of the definition of $f, g$ and $\mathcal{D}(u,v)$ given in \eqref{fg} and \eqref{D}, respectively, we have exploited that $\mathcal{D}(u,v)=\chi \|f\|_{L^2(\Omega)}^2 + \|g\|_{L^2(\Omega)}^2$.
\end{proof}
\end{lemma}
Now, by applying Lemma \ref{DerivativeEstimate} and Lemma \ref{TheoremBU} we derive a differential inequality for $-\mathcal{F}_h(u(\cdot,t), v(\cdot,t))$ with superlinearly growing nonlinearity. This implies that $(u,v)$ cannot exist globally for initial data $(u_0,v_0)$ with large negative energy $\mathcal{F}_h(u_0,v_0)$.
\begin{lemma}\label{LemmaBU}
Let $m, A, h, \chi>0$, $\kappa > n-2$ if $n\geq 3$, and $\kappa>2$ if $n=2$. Then there exist $K=K(m,A,\kappa, h, \chi)>0$ and $C=C(m,A,\kappa, h, \chi)>0$ such that for each $(u_0,v_0)$ from the set
\begin{equation*} 
\begin{split}
\mathcal{\tilde{B}}(m,A,\kappa, h, \chi):= &\left\{(u_0, v_0) \textrm{ complying with } \eqref{AssumpAnnex}, \textrm{ and with } \int_\Omega u_0=m, \|v_0\|_{W^{1,2}(\Omega)} \leq A  \textrm{ and } \mathcal{F}_h(u_0,v_0) \leq - K(m, A, \kappa, h, \chi)\right\},
\end{split}
\end{equation*}
the corresponding solution $(u,v)$ of \eqref{problemB} has the property 
\begin{equation}\label{ODI}
\mathcal{F}_h(u(\cdot,t), v(\cdot,t)) \leq \frac{\mathcal{F}_h(u_0, v_0)}{(1-Ct)^{\frac{\theta}{1-\theta}}} \quad \textrm{for all } t \in (0, \TM),
\end{equation}
where $\theta \in \left(\frac{1}{2},1\right)$ is as given by \eqref{theta}.
In particular, for any such solution we have $\TM<\infty$, i.e., $(u,v)$ blows up in finite time.
\begin{proof}
By reasoning as in \cite[Lemma 5.2]{W}, exploiting Lemma \ref{UpperBoundV} and bounds \eqref{MassZero} and \eqref{VZero}, we have that the solution $(u,v)$ of problem \eqref{problemB} emanating from $(u_0,v_0)$ belongs to the set $\mathcal{S}(m,M,B, \kappa, h)$ for all $t \in (0,\TM)$. This implies that Lemma \ref{TheoremBU} guarantees the existence of  $\const{Ab1}=\const{Ab1}(m,M,B, \kappa, h,\chi)$ such that the solution $(u,v)$ fulfills 
\begin{equation}\label{EstFD1}
\mathcal{F}_h(u,v) \geq -\const{Ab1} \, (\mathcal{D}^{\theta}(u,v)+1) \quad \textrm{for all } t \in (0,\TM).
\end{equation}
To show that \eqref{ODI} holds for all $(u_0,v_0) \in \mathcal{\tilde{B}}(m,A,\kappa, h, \chi)$, let us define 
\begin{equation}\label{KC}
K=K(m,A,\kappa, h, \chi):= 2 \const{Ab1} \quad \textrm{and} \quad C=C(m,A,\kappa, h, \chi):= \frac{1-\theta}{2 \const{Ab1} \theta}.
\end{equation}
From the regularity properties of $(u,v)$ 
we derive that 
\begin{equation*}
y_h(t):= - \mathcal{F}_h(u(\cdot,t), v(\cdot,t))\in C([0,\TM))\cap C^1((0,\TM)).
\end{equation*}
Moreover $y_h(t)$ is positive; indeed from
relation \eqref{DerEst1} we have that $y_h$ is nondecreasing and since $y_h(0)\geq K=2 \const{Ab1}$ (recall \eqref{KC}), this leads to
\begin{equation}\label{yPos}
y_h(t) \geq y_h(0) \geq 2 \const{Ab1}>0 \quad \textrm{for all } t \in (0,\TM).
\end{equation}
By manipulating inequality \eqref{EstFD1} and taking into account bound \eqref{yPos}, we obtain
\begin{equation*}
\mathcal{D}^{\theta}(u(\cdot,t), v(\cdot,t)) \geq \frac{y_h(t)}{\const{Ab1}}-1\geq \frac{y_h(t)}{2\const{Ab1}} \quad \textrm{for all } t \in (0,\TM),
\end{equation*}
that, inserted in relation \eqref{DerEst1}, implies
\begin{equation*}
y'_h(t) \geq \left(\frac{y_h(t)}{2\const{Ab1}}\right)^{\frac{1}{\theta}} \quad \textrm{for all } t \in (0,\TM).
\end{equation*}
After some computations we have
\begin{equation*}
y_h(t) \geq y_h(0) \left(1-\frac{1-\theta}{\theta}(2\const{Ab1})^{-\frac{1}{\theta}}  y_h^{\frac{1-\theta}{\theta}}(0) t\right)^{\frac{-\theta}{1-\theta}} \quad \textrm{for all } t \in (0,\TM),
\end{equation*}
which gives inequality \eqref{ODI} after have noted from relation \eqref{yPos} (and recalling definition \eqref{KC}) that 
\begin{equation*}
\frac{1-\theta}{\theta}(2\const{Ab1})^{-\frac{1}{\theta}} y_h^{\frac{1-\theta}{\theta}}(0) \geq \frac{1-\theta}{2 \const{Ab1}\theta}=C, 
\end{equation*}
so we can conclude.
\end{proof}
\end{lemma}
Finally, we can state the other main result of this Appendix:
\begin{theorem}[Finite-time blow-up $\tau=1$]\label{BlowUpParabolic}
Let $\Omega=B_R \subset \mathbb{R}^n$, with $n \geq 2$ and $R>0$, and let $m, A, h, \chi>0$. Then there exist $T(m, A, h, \chi)>0$ and $K(m, A, h, \chi)>0$ with the property that given any $(u_0,v_0)$ from the set
\begin{equation*} 
\begin{split}
\mathcal{B}(m,A,h,\chi)= &\left\{(u_0, v_0) \textrm{ complying with } \eqref{AssumpAnnex}, \textrm{ and with }  \int_\Omega u_0=m, \|v_0\|_{W^{1,2}(\Omega)} \leq A  \textrm{ and } \mathcal{F}_h(u_0,v_0) \leq - K(m,A,h, \chi)\right\},
\end{split}
\end{equation*}
for the corresponding solution $(u,v)$ of \eqref{problemB} we have $\TM \leq T(m, A, h, \chi) < \infty$, i.e., $(u,v)$ blows up before or at time $T(m,A, h, \chi)$.
\begin{proof}
In order to remove the auxiliary parameter $\kappa$, that appears in Lemma \ref{LemmaBU}, we have to fix an arbitrary $\kappa>n-2$ if $n\geq 3$ and a $\kappa>2$ if $n=2$. Therefore, we apply Lemma \ref{LemmaBU} to derive the statement for $K(m,A,h, \chi):=K(m,A, \kappa, h, \chi)$ and $T(m,A,h, \chi):=\frac{1}{C(m,A, \kappa, h, \chi)}$ with $K(m,A, \kappa, h, \chi)$ and $C(m,A, \kappa, h, \chi)$ as given by Lemma \ref{LemmaBU}.  
\end{proof}
\end{theorem}
\begin{remark}[Comparison of the time $T$ of blow-up: Neumann vs Robin]
Let us assume for simplicity $\chi=1$ and denote by $T_N$ and $T_R$ the time of blow-up in the Neumann and in the Robin case, respectively. By \cite[Theorem 1.1]{W} we know that 
\[
T_N= \frac{1}{C(m,A, \kappa)}=\frac{2c_N \theta}{1-\theta},
\] 
where $\theta$ is given in \eqref{theta} and $c_N:=\max\{\overline{C}(m,A,\kappa), \underline{C}(m,A,\kappa)\}$. 
By retracing the previous lemmas from Lemma \ref{FirstEstimateInt_uv} and 
paying attention to the dependence of the constants on $h$ (see in particular estimate \eqref{FirstEst3}, we derive that 
\[
T_R= \frac{1}{C(m,A, \kappa,h)}=\frac{2c_R \theta}{1-\theta},
\] 
where $c_R:=\max\{\overline{C}(m,A,\kappa), \underline{C}(m,A,\kappa)
\rho(h)\}$ with $\rho(h):=\left(1+\frac{hn}{R}+3h^2\right)^{\frac{n+2}{2}}$.
By comparing the values $T_N$ and $T_R$ we obtain that 
\[
c_N:=\max\{\overline{C}(m,A,\kappa), \underline{C}(m,A,\kappa)\}
< c_R:=\max\{\overline{C}(m,A,\kappa), \underline{C}(m,A,\kappa)\rho(h)\},
\]
this implies that 
\[
T_N < T_R,
\]
and they coincide when $h \to 0$ since $\rho(h) \to 1$.
\end{remark}

\subsubsection*{Acknowledgements}
The authors SF and GV are members of the Gruppo Nazionale per l'Analisi Matematica, la Probabilit\`{a} e le loro Applicazioni (GNAMPA) of the Istituto Nazionale di Alta Matematica (INdAM) and they are partly supported 
by MIUR (Italian Ministry of Education, University and Research) Prin 2022 \textit{Nonlinear differential problems with applications to real phenomena} (Grant Number: 2022ZXZTN2). 
YT is partially supported by JSPS KAKENHI (Grant Number: JP24K22844). 
GV is also supported by \textit{Partial Differential Equations and their role in understanding natural phenomena} (2023), funded by \href{https://www.fondazionedisardegna.it/}{Fondazione di Sardegna} (CUP F23C25000080007). KB is partially supported by Pasargad Institute for Advanced Innovative Solutions (PIAIS, Grant Number 1403-10193).

\end{document}